\documentclass[10pt]{article}
\usepackage{amsmath}
\usepackage{amsthm}
\usepackage{graphicx}
\usepackage{amsfonts}
\usepackage{amssymb}
\usepackage{enumerate}
\usepackage[margin=0.85in]{geometry}
\usepackage{comment}
\usepackage{xcolor}
\usepackage{hyperref}
\usepackage{cleveref}
\usepackage{upgreek}
\usepackage{todonotes}
\usepackage{tikz}
\usepackage{tikz-cd}
\usetikzlibrary{matrix,arrows,decorations.pathmorphing, shapes, positioning,calc}
\usepackage{graphicx}
\usepackage{float}
\usepackage[export]{adjustbox}

\usepackage[ruled]{algorithm2e}

\usepackage[english]{babel}
\usepackage[latin1]{inputenc}
\usepackage[numbers]{natbib}

\usepackage[normalem]{ulem}
\usepackage{bbm}

\DeclareOldFontCommand{\bf}{\normalfont\bfseries}{\mathbf}
\DeclareOldFontCommand{\cal}{\normalfont\bfseries}{\mathcal}

\usepackage{xcolor}
\pagecolor{white}

\newtheorem{theorem}{Theorem}[section]

\newtheorem{proposition}[theorem]{Proposition}

\theoremstyle{definition}

\newtheorem{definition}[theorem]{Definition}
\newtheorem{remark}[theorem]{Remark}
\newtheorem{condition}[theorem]{Condition}
\crefname{condition}{condition}{conditions}
\numberwithin{equation}{section}

\newcommand{\diff}{\,\mathrm{d}}

\newcommand{\ceil}[1]{\lceil{#1}\rceil}
\renewcommand{\d}{\mathrm{d}}
\newcommand{\e}{\mathrm{e}}

\newcommand{\opt}{\widetilde{\pi}}
\newcommand{\paran}[1]{\left({#1}
\right)}
\newcommand{\abs}[1]{\lvert#1\rvert}

\newcommand\cA{\mathcal A}
\newcommand\cB{\mathcal B}
\newcommand\cC{\mathcal C}

\newcommand\cE{\mathcal E}
\newcommand\cF{\mathcal F}

\newcommand\cI{\mathcal I}

\newcommand\cN{\mathcal N}

\newcommand\EE{\mathbb E}

\newcommand\RR{\mathbb R}

\def \E{\mathbb{E}}
\def \F{\mathbb{F}}

\def \LL{\mathbb{L}}

\def \N{\mathbb{N}}

\def \P{\mathbb{P}}
\def \Q{\mathbb{Q}}
\def \R{\mathbb{R}}

\definecolor{antiquefuchsia}{rgb}{0.57, 0.36, 0.51}
\definecolor{airforceblue}{rgb}{0.36, 0.54, 0.66}
\definecolor{asparagus}{rgb}{0.53, 0.66, 0.42}
\definecolor{blond}{rgb}{0.98, 0.94, 0.75}
\definecolor{amber(sae/ece)}{rgb}{1.0, 0.49, 0.0}
\definecolor{cadmiumorange}{rgb}{0.93, 0.53, 0.18}
\definecolor{denim}{rgb}{0.08, 0.38, 0.74}
\definecolor{amber}{rgb}{1.0, 0.75, 0.0}
\definecolor{chromeyellow}{rgb}{1.0, 0.65, 0.0}
\definecolor{cobalt}{rgb}{0.0, 0.28, 0.67}

\pdfoutput=1
\usepackage{color}
\usepackage{subfigure}
\begin{document}

\title{A Deep Learning Method for Optimal Investment Under Relative Performance Criteria Among Heterogeneous Agents
}
\author{Mathieu Lauri\`ere\footnote{Shanghai Frontiers Science Center of Artificial Intelligence and Deep Learning; NYU-ECNU Institute of Mathematical Sciences, NYU Shanghai, 567 West Yangsi Road, Shanghai, 200126, People's Republic of China, mathieu.lauriere@nyu.edu} $\qquad$ Ludovic Tangpi\footnote{Princeton University, ORFE, PACM and BCF. Princeton, 08544, NJ. USA, ludovic.tangpi@princeton.edu. Funding: NSF grant DMS- 2005832, NSF CAREER award DMS-2143861 and AMS Claytor-Gilmer fellowship.}
$ \qquad$  Xuchen Zhou \footnote{Princeton University, ORFE. Princeton, 08544, NJ. USA, xuchen.zhou@princeton.edu. Funding: NSF grant DMS- 2005832.}}

\maketitle

\paragraph{\textbf{Abstract}} 
Graphon games have been introduced to study games with many players who interact through a weighted graph of interaction. By passing to the limit, a game with a continuum of players is obtained, in which the interactions are through a graphon. In this paper, we focus on a graphon game for optimal investment under relative performance criteria, and we propose a deep learning method. The method builds upon two key ingredients: first, a characterization of Nash equilibria by forward-backward stochastic differential equations and, second, recent advances of machine learning algorithms for stochastic differential games. We provide numerical experiments on two different financial models. In each model, we compare the effect of several graphons, which correspond to different structures of interactions.



\section{Introduction}

Optimal investment constitutes a pivotal facet of the financial industry. 
As market dynamics grow progressively complex, the imperative to develop advanced and sophisticated models becomes more pronounced. Equally crucial is the development of efficient numerical simulation algorithms for optimal portfolio management. 
In fact hedge funds invest in diverse financial products, impose numerous constraints on their portfolios, and necessitate frequent update of their investment strategies.

In the modern financial landscape, the intricacy of simulating optimal portfolios is further compounded by the imperative for hedge funds to compete vigorously with their counterparts. Hedge funds usually raise capital from the same pool of investors, compelling them to consistently outperform their peers to attract and retain clients. 
Mathematically, such portfolio optimization problem can be stated as a stochastic differential game.
This is a very difficult problem in practice particularly given the exponential growth in the number of hedge funds in recent years.

One of the first papers treating the problem in a systematic, rigorous way is the work of \citet{espinosa_touzi_2013} on exponential utility maximization under relative performance concern, who showed that deriving Nash equilibria for such a game reduces to solving an $n$-dimensional system of backward stochastic differential equations (BSDEs), where $n$ is the number of competing firms.
The analysis of such equations has given rise to a sizable literature due to the complexity and relevance of the problem. See for instance \cite{luo2019multidimensional,frei2014splitting,Frei-Reis11,Xing-Zit18}.
It is well-known that the computation of solutions for large dimensional BSDEs is usually a complex problem, and very few works have attempted to solve such high-dimensional BSDEs in the context of games using traditional numerical methods such as least-square regression, quantization or tree-based approximations. See e.g.~\cite{chessari2023numerical} for an overview of numerical methods for BSDEs.

A key insight is to consider the limiting case of a game with infinitely many firms.
For perfectly symmetric games where firms are (probabilistically) identical, the $N=\infty$ regime is known as mean field games as introduced in \cite{MR2295621,huang2007invariance}.
This idea was first considered in \cite{lacker2019mean}.
In this case, the problem of finding equilibria reduces to solving a one-dimensional equation.
We refer to \cite{Lacker-Soret20,hu2022n,dos2019forward,dos2022forward,guanxing2020mean} for a small sample of works on the problem.
In particular, the numerical simulation of mean field game equilibria becomes very promising as the issue of dimensions is no-longer present.
Moreover, and crucially, it is well-known that mean field games equilibria allow to obtain $\varepsilon$-Nash equilibria\footnote{Roughly speaking, a family of strategies is an $\varepsilon$-Nash equilibrium if single a player deviating can at most increase their reward by $\varepsilon$.} for the finite player game with a large enough population.

The goal of the present paper is to propose an algorithm to efficiently simulate mean field equilibria in a general utility maximization problem with relative performance concern.
In order to make the game considered here more realistic, we will further consider the case where participating firms are not completely symmetric.
In fact, we will allow each firm to compete with only a smaller group of other market participants.
This reflects for instance the fact that hedge funds compete mostly with those investing in similar asset classes or with funds of similar size.
Such games are known as games on graphs and, as showed by \citet{tangpi-zhou22}, in this case the game in the $N=\infty$ regime is no longer a mean field game, but a so-called graphon game; see Section \ref{sec: graphon game} for definition.
Essentially, graphons are symmetric, measurable functions from $[0,1]\times[0,1]$ to $[0,1]$ representing natural limits of (dense) graphs, see \cite{Lovasz}.
It is therefore natural that such object will capture the interaction among players in the limiting regime of games on graphs.
Note that mean field games are simple examples of graphon game.
We will thus propose a numerical simulation method for this more general type of problems.
Papers dealing with stochastic differential games in a rather generic context include \cite{caines2019graphon,aurell2022finite,aurell2022stochastic,Lacker-Soret22}.
These papers are mostly concerned with existence and uniqueness of equilibria.
The paper of \citet{tangpi-zhou22} presents the problem of exponential utility maximization among heterogenous, competing agents as a graphon game.
We will make use of the setting and results of the latter reference, and focus on the issue of efficient simulation of equilibria in the general context.
Despite a huge effort to understand simulation of stochastic differential games, the issue of simulating graphon games has been very sparsely discussed in the literature, with the exception of~\cite{aurell2022finite,cuilearning,Fab-Cui-Koe23}.

The computational method we use builds upon the fact that graphon equilibria are fully characterized by solutions of graphon BSDE of McKean-Vlasov type \cite{aurell2022stochastic,tangpi-zhou22,Bayr-Wu-Zhang22}.
Graphon BSDEs are systems of BSDEs indexed by $u \in [0,1]$, that are fully coupled (through a graphon) and involve the law of the solution.
In particular, we will need to solve a continuum of coupled McKean-Vlasov BSDEs. 
To this end, we will learn the solution as a parameterized function of the index, leveraging recent machine learning methods for control problems and games (see e.g.~\cite{hu2023recent} for an overview). More precisely, we use a shooting method and rewrite the FBSDE as an optimal control problem for two forward SDEs with a terminal cost encouraging the second process, $Y$, to satisfy the correct terminal condition. In this new problem, the controls are the initial value of $Y_0^u = y_0(u,X_0^u)$ and the volatility $Z_t^u = z(t,u,X_t^u)$ for the second process. The graphon interactions are replaced by Monte Carlo samples. The shooting method for BSDEs using neural networks was proposed in~\cite{han2018solving} to approximate the solution to PDEs. It was then extended in~\cite{carmona2022convergence} to FBSDEs of MKV type. A similar shooting method was used in~\cite{aurell2022finite} to solve graphon forward-backward ODEs for finite-state games, by approximating the distribution and the value function with neural networks taking the player's index as an input. The idea of learning parameterized functions of the index has then been used in~\cite{cuilearning,Fab-Cui-Koe23} but still for finite-state graphon games. Here, we propose an algorithm to solve continuous space graphon games using a combination of the aforementioned methods (FBSDEs, shooting method and neural network taking the player's index as an input).

The method just described allows to efficiently simulate the equilibrium strategies along with the equilibrium wealth dynamics and optimal utilities.
We do so for various interaction graphons and financial models.
Since in the case of deterministic coefficients we have equilibria in closed form, we can perform a sanity check of our method and compute the estimation errors along with the estimated losses.
In particular, for the case where the graphon game reduces to a mean field game, the validation loss is of the order $10^{-11}\%$ and the relative approximation error is of the order $10^{-7}\%$.
Our simulations give interesting insights in the utility maximization game.
We present how the kind of interaction among the players impacts the investment strategy.
For instance, it is interesting to notice that players facing interaction with a larger population (i.e. facing more competition pressure) have lower optimal utility.

In the rest of the paper, we begin by clearly presenting the model as well as the optimization game we are interested in.
We also recall some of the theoretical characterization results that we use in our algorithm.
In the following section we describe the simulation methods and the neural network architecture that is used.
In the final section of the paper we present simulation results for various models and graphons.
We also discuss the financial insights gained from our experiments.
 
\section{Utility maximization under relative performance concern}

The underlying motivation of this paper is to present a simulation method for the equilibrium strategy as well as the equilibrium wealth of a large population of competitive, heterogeneous investors.
This will be done using the recent notion of graphon games, see theoretical results derived by \citeauthor*{tangpi-zhou22} \cite{tangpi-zhou22}, as well as a neural network approach to solving McKean-Vlasov BSDEs.
Before introducing the graphon game per se, let us discuss its finite population counterpart.
This is the actual problem being modeled.
We will later explain how it relates to the graphon game that we will simulate.

Let $d \in \N^\star$ be fixed and equip the canonical space $\cC_d$ of $\R^d$--valued continuous functions on $[0,T]$ with the Wiener measure $\P^o$ and canonical process $W$.
We fix as index space the unit interval $I:=[0,1]$ and denote by $(I, \mathcal{B}(I), \lambda_I)$ the Lebesgue space.
Let $I =[0,1]$ be fixed and equip the unit interval $I$ with the Borel $\sigma$--algebra $\mathcal{B}(I)$ and Lebesgue measure $\mu_I$.
By \cite[Theorem 1]{sun2006exact}, there is a probability space $(I, \mathcal{I}, \mu)$ extending the Lebesgue space $(I, \mathcal{B}(I), \mu_I)$, a probability space $(\Omega ,{\mathcal{F}}, \P)$ and a \emph{Fubini extension} $(\Omega\times I, {\mathcal{F}}\boxtimes \mathcal{I}, \P\boxtimes \lambda)$ of the product space $(\Omega\times I, {\mathcal{F}}\otimes \mathcal{I}, \P\otimes \lambda)$ on which can be constructed an ${\mathcal{F}}\boxtimes \mathcal{I}$--measurable process $\mathbb{W}:\Omega \times I \longrightarrow \cC_{d}$ with $\mathbb{W}(\cdot, u) = W^u(\cdot)$, such that the family $(W^u)_{u \in I}$ is \emph{$\P$--essentially pairwise independent} ($\P$--e.p.i. for short)\footnote{A family of random variables $(\zeta^u)_{u\in I}$ is $\tilde\P$--e.p.i. if for $\lambda$--a.e. $u \in I$, $\zeta^x$ is $\tilde\P$-independent of $\zeta^y$ for $\lambda$--a.s. $y \in I$.} and with law $\P^o$ for all $u \in I$.
In particular, $(W^u)_{u\in I}$ is a family of $\P$--e.p.i Brownian motions on $(\Omega, \cF,\P)$.
We refer the unfamiliar readers to \citeauthor*{sun2006exact} \cite{sun2006exact} for an introduction to Fubini extensions.

\subsection{Market model and the finite population game  \label{sec: n-agent game}}
Let $N\in \mathbb{N}^\star$ be fixed and let $(u^i)_{i\in \{1,\dots,N\}}$ be a given sequence of elements of the unit interval $I$. Assume for simplicity that $(u^i)_{i\in\{1,\dots,N\}}$ are such that $(W^{u^{1}},\dots,W^{u^N})$ denotes a family of $\P$-independent $\P$--Brownian motions.
Note that when there is no ambiguity ,we will index players in the $N$--player games by $\{1,\dots,N\}$ rather than $\{u^1,\dots,u^N\}$, and thus write
\begin{equation*}
    W^i \equiv W^{u^i}.
\end{equation*}

We will denote by $\F^i$ the completed filtration generated by $W^i$, and $\F_N$ the completed filtration generated by $(W^1,\dots, W^N)$.
The financial market consists of $N$ agents trading in a common risk--less bond with interest rate $r=0$ and $N\times d$ stocks.
That is, agent $i$ trades in $d$ stocks whose prices are represented by a vector $S^i \in \RR^d$ following the dynamics
\begin{align}\label{stock evolution}
    \d S_t^{i} = \mathrm{diag}(S^i_t)\big(\mu^{i}_t \d t + \sigma_t^i\d W^{i}_t \big) \quad i=1,\dots,N,
\end{align}
where $\mathrm{diag}(x)$ is the diagonal matrix with entries $x\in \R^d$ on the diagonal, and where $\mu^i$ and $ \sigma^i$ are sufficiently integrable, $\F^i$--adapted processes taking values in $\R^d$ and $\R^{d\times d}$, respectively.
Given a trading strategy $\pi$, which we assume to be a $\F_N$--predictable and $\R^d$--valued process representing the amount of money invested in the stock market, we denote by $X_t^{i, \pi}$ the wealth of agent $i$ at time $t$.
When starting with the initial position $\xi^i$ and assuming $\pi$ to be self--financing, $X^{i,\pi}$ satisfies 
\begin{align*}
     \d X_t^{i, \pi} =   \pi_t \cdot \left(\sigma_t^{i} \theta_t^i \d t + \sigma_t^i\d W_t^i   \right),\quad X^{i,\pi}_0 = \xi^i
\end{align*}
where we introduced the process $\theta^i$ given by 
$$\theta_t := (\sigma_t^i)^{-1}\mu_t^i.
$$ 
The optimization problem that a generic agent $i$ faces is 
\begin{align}\begin{split}
\label{eq: n-agent obj}
 V^{i,N}
    &:=V^{i,N}((\pi^j)_{j\neq i})  
    := \sup_{\pi \in \cA_i} \E \bigg[ -\exp\bigg\{-\frac{1}{\eta^i}\bigg(X_T^{i,\pi} - \frac{\rho}{N-1}\sum_{j\in \{1,\dots,N\}\setminus \{i\}}\lambda_{ij} X_T^{j,\pi^j}\bigg)\bigg\} \bigg]
    \end{split}
\end{align}
for a risk aversion parameter $\eta_i > 0$, a matrix $(\lambda_{ij}^N)_{(i,j)\in \{1,\dots,N\}^2}$ and a fixed parameter $\rho\in [0,1]$.
We will define the set of admissible controls $\mathcal A_i$ below.
Problem \ref{eq: n-agent obj} is the utility maximization of an agent who is concerned with the relative performance of \emph{some of} their peers.
The underlying idea here is that while financial investors are in non--cooperative competition, they do not necessarily compete with all market participants.
For instance, most often, hedge funds will compete only with competitors raising capital from similar clients or those investing in the same sectors.
The parameter $\rho$ (which we will often take to be $1$) models how much market participants prioritize competition, whereas $\lambda_{ij}\in \{0,1\}$ indicates whether agent $i$ is concerned with the performance of agent $j$.
The matrix $(\lambda_{ij})_{(i,j)\in \{1,\dots,N\}^2}$ should be thought of as the adjacency matrix of the graph on which the $N$--player game is written.
The set $\cA_i$ of admissible trading strategies for agent $i$ will be defined as follows.
\begin{definition}[Admissibility]
\label{defn: admissible strategy}
    Let $A_i$ be a closed convex subset of $\R^d$.
    A strategy $\pi^i$ is admissible if $\pi^i$ is $\F_N$--progressive, $A_i$--valued and square integrable, and furthermore, for every $j\in \{1,\dots,N\}$ there is $p>2$ such that the family
    \[
        \big\{\e^{\frac{p}{\eta^i}\rho \lambda_{ij}X_{\tau}^{i,\pi^i}};\; \text{with $\tau$ a $\F_N$--stopping time on $[0,T]$}  \big\}
    \]
    is uniformly integrable. We denote by $\cA_i$ the set of admissible strategies for player $i$. 
\end{definition}
When  $\rho = 0$, this problem reduces to the standard utility maximization problem studied for instance by \cite{hu2005utility,Rouge-Karoui01,tarpooptim,Delbaen-expo-uti} and many others.
With $\lambda_{ij}=1$ for all $i,j\in \{1,\dots,N\}$, this is the (standard) utility maximization problem under relative performance concern first studied by \cite{EspinosaThesis,espinosa_touzi_2013} and later extended by \cite{hu2022n,anthropelos2022competition,Lacker-Soret20,lacker2019mean,guanxing2020mean}.

We will be interested in computing Nash equilibria for the game described above: that is, strategies $(\widetilde \pi^1,\dots, \widetilde \pi^N)$ such that for each $i$
\begin{equation*}
    V^{i,N}((\widetilde\pi^j)_{j\neq i}) = \E \bigg[ -\exp\bigg\{-\frac{1}{\eta^i}\bigg(X_T^{i,\widetilde\pi^i} - \frac{\rho}{N-1}\sum_{j\in \{1,\dots,N\}\setminus \{i\}}\lambda_{ij} X_T^{j,\widetilde\pi^j}\bigg)\bigg\} \bigg].
\end{equation*}
 Existence of Nash equilibria of this non--cooperative game (at least when $\lambda^N_{ij}=1$ for all $(i,j)\in \{1,\dots,N\}^2$) is a well--studied topic, refer again to \cite{EspinosaThesis,espinosa_touzi_2013, hu2022n,anthropelos2022competition,Lacker-Soret20,lacker2019mean,guanxing2020mean}.
 However, to the best of our knowledge, the issue of numerical computation of a Nash equilibrium for such utility maximization problems has not been addressed despite it central importance for practical applications.
Typically, the computational cost of computing an exact Nash equilibrium becomes prohibitive when the size of the population is very large, as is required in applications.
For instance, we think of $N$ here as the number of hedge funds operating on the U.S. market.
Many financial publications estimate there to be more than 3.500 hedge funds operating in the U.S. as of 2023.

These computational challenges have motivated the development of approximation methods that can provide a tractable framework while providing good approximate solutions. 
Among the methods that have been used to get around such issues in stochastic differential games, a powerful one consists in considering the $N\to\infty$ version of the game known as the \emph{mean field game}.
However, mean field games are limits of \emph{symmetric} games, in which each player interacts only with the empirical distribution of agents. This is not the case for the game considered in the present work due to the terms $(\lambda_{ij})_{i,j}$.
The natural limit of the game introduced above with heterogeneous interactions is a graphon game.

\subsection{The graphon game \label{sec: graphon game}}
The graphon game is a game with a continuum of players.
Players are indexed by $u \in I$ and \emph{interaction} among the continuum of agents will be modeled by a graphon, which is a symmetric and measurable function
\begin{equation*}
    G:I\times I\to I.
 \end{equation*} 
 In general, graphons can take any real values but in this work we focus on $[0,1]$-valued graphons. 
Throughout the paper, we fix a graphon $G$.
Concrete examples will be studied in the next sections.
Assume that each player $u$ invests in a stock $S^u$ with dynamics
\begin{align*}
    \d S_t^{u} = \mathrm{diag}(S^u_t)\big(\mu^{u}_t \d t + \sigma_t^u\d W^{u}_t\big), \quad S_0^u > 0,
\end{align*}
so that when agent $u$ employs a strategy $\pi$, her wealth follows the dynamics
\begin{align}
\label{eqn:mf-portfolio}
   \d X_t^u &= \pi_t \cdot \left(\sigma^u_t\theta^u_t \d t  + \sigma_t^u\d W_t^u 
   \right), \quad X^u_0 = \xi^u
\end{align}
where $\sigma$, $\theta$, are $\mathcal{B}([0,T])\times\cI \boxtimes \cF$--measurable stochastic processes and such that $(\theta^u)_{u\in I}$ and $(\sigma^u)_{u\in I}$ are e.p.i. and identically distributed. We will denote by $\F^u$ the completed filtration generated by $W^u$, and impose the following restrictions on the initial wealth and stock price coefficients:
\begin{condition} \label{cdn:type vector} ~\
\begin{enumerate}
    \item The function $(u,\omega)\mapsto \xi^u(\omega)$ is measurable and $\E[|\xi^u|^4]<\infty$ for all $u\in I$.
    \item $\mu^u,\sigma^u$, and $\theta^u$ are bounded uniformly in $u \in I$.
    \item For $\mu$-almost every $u \in I$, $\mu^u$, $\sigma^u$ and $\theta^u$ are $\F^u$-predictable.
    \item For every fixed $(t,u) \in [0,T]\times I$, $\sigma_t^u$ is uniformly elliptic, that is, $\sigma_t^u(\sigma_t^u)^\top\geq c I_{d\times d}$ in the sense of positive semi-definite matrices, for some $c >0$ independent of $(t,u)$.
\end{enumerate}
\end{condition}

Let a measurable mapping $I\ni u \mapsto \eta^u\in (0,\infty)$ define a family of risk-aversion parameters for the players.
Admissible trading strategies in the graphon game are defined as follows:
\begin{definition}
    \begin{itemize}
        \item A strategy profile is a family $(\pi^u)_{u\in I}$ of processes taking values in $\R^d$ such that for every $u$, $\pi^u$ is $\F^u$--progressive, $(t,u,\omega)\mapsto \pi^u_t(\omega)$ is $\cB([0,T])\otimes\cI\boxtimes \cF$--measurable. 
        \item For each $u \in I$ let $A_u$ be a closed subset of $\R^d$.
        A  strategy $\pi^u$ is admissible if it is $\F^u$--progressive and takes values in $A_u$; and there is a strategy profile $(\pi^u)_u$ such that $\int_I\E[\int_0^T\|\pi^u_t\|^2\diff t]\mu(\diff u)<\infty$.        
        We denote by $\pi^u\in \cA^G_u$ the set of admissible strategies for player $u$ in the graphon game.
    \end{itemize}
\end{definition}
Notice that the definition of admissible strategies in the graphon game is simpler than the analogous definition in the finite player game, see Definition~\ref{defn: admissible strategy}. Indeed, in the graphon game, the player's states become e.p.i. so the strategy of player $u$ is required to be progressively measurable simply with respect to its individual filtration $\F^u$. Although this is only a small simplification from the mathematical viewpoint, it has its importance for the design of numerical methods since it implies that it is sufficient to focus on strategies that depend only on the individual player's state. 
The graphon utility maximization problem for player $u \in I$ now takes the form
\begin{align}
    \begin{split}
    \label{eq:graphon-obj}
    V^{u,G} &= V^{u,G}\left(\{ \pi^v \}_{v \neq u}\right) 
    := \sup_{\pi^u \in \cA^{G}_u} \E\left[-\exp\left(-\frac{1}{\eta^u}\left(X_T^{u,\pi^u} - \rho \E\big[ \int_I  X_T^{v, {\pi}^v} G(u,v) \d v \big] \right)\right) \right].
\end{split}\end{align}
In particular, we will be interested in the graphon equilibria:
\begin{definition}
    A family of admissible strategy profiles $(\widetilde \pi^u)_{u\in I}$ is called a graphon (Nash) equilibrium if {}for almost every $u\in I$, 
    it holds
        \begin{equation*}
            V^{u,G}\left(\{ \widetilde \pi^u \}_{v \neq u}\right) \\
     :=  \E\left[-\exp\left(-\frac{1}{\eta^u}\left(X_T^{u,\widetilde\pi^u} - \rho \E\big[ \int_I  X_T^{v, {\widetilde\pi}^v} G(u,v) \d v  \big] \right)\right) \right].
        \end{equation*}
\end{definition}

\subsection{Characterization, convergence and approximation}
As the size of the population grows larger, the finite agent game described in \Cref{sec: n-agent game} converges to the graphon game described in \Cref{sec: graphon game} in terms of both equilibrium strategies and equilibrium utilities. This result will make the numerical approximation of the finite-agent game possible since, given that population is large, we can now attempt to numerically solve for the forward-backward stochastic differential equations (FBSDEs) characterizing the graphon equilibrium, which is much more tractable. Such convergence result is guaranteed under some regularity constraints regarding the stock price coefficients and the graphon. In the rest of this section, we will first present \Cref{prop:graphon-bsde}, which allows us to characterize the graphon Nash equilibrium in terms of a graphon FBSDE. This is a system of (infinitely many) coupled FBSDEs. The graphon FBSDEs in \Cref{prop:graphon-bsde} will serve as a starting point for our approximation method. We will then conclude this section by presenting the convergence result from the $n$-agent game to the graphon game to justify our choice of computing the graphon equilibrium. 

For computation purposes, the following condition will be held true for the rest of this work
\begin{condition}\label{cdn: unconstrained}
The strategies are unconstrained, i.e., $A^u = \R^d$ for all $u \in I$.
\end{condition}

\begin{proposition}
\label{prop:graphon-bsde}
   Let \Cref{cdn:type vector} and \Cref{cdn: unconstrained} be satisfied. The graphon game described in~\eqref{eq:graphon-obj} admits a graphon Nash equilibrium $(\tilde \pi^u)_{u\in I}$ such that for almost every $u \in I$ it holds  
    \begin{align}
    \label{eq:optim.sol.no.common}
        \opt^u_t = (\sigma_t^{u})^{-1}\paran{Z_t^{u}+\eta^u\theta_t^u}\,\, \mathrm{d}t\otimes\d u\boxtimes\P\text{--a.s.} \text{ and}\quad V^{u,G} = -\exp\Big(- \frac{1}{\eta^u} \Big(\xi^u - \int_I \E[\rho \xi^v]G(u,v) \d v - Y_0^u \Big) \Big)
    \end{align}
    where $(X^u, Y^u, Z^u)_{u\in I}$ solves the following FBSDE system
    \begin{align}   \label{eqn:BSDE-graphon.decoupled}
    \begin{cases} 
        \d X_t^u = \widetilde \pi_t^u \cdot (\sigma_t^u\theta_t^u \d t + \sigma_t^u \d W_t^u), \quad \widetilde \pi_t^u = (\sigma_t^u)^{-1}\big(Z_t^u+ \eta^u\theta_t^u\big), \,\quad \d t\otimes\d u\boxtimes\P\text{--a.s.} \\ 
        \d Y_t^u =  \bigg(\frac{\eta^u}{2}\abs{\theta_t^u}^2 + Z_t^{u}\cdot \theta_t^u - \E\Big[\int_I \rho \paran{Z_t^{v}+\eta^v\theta_t^v}\cdot \theta^v_tG(u,v)\d v\Big] \bigg) \d t \\
        \qquad \quad+  Z_t^{u} \cdot \d W_t^u, \, \quad \d t\otimes \d u\boxtimes \P\text{--a.s.,}    \\
        Y_T^u = 0, \quad X_0^u = \xi^u.
    \end{cases}
    \end{align} with $(u,t,\omega)\mapsto Z^u_t(\omega)$ measurable and $(X^u,Y^u, Z^u)$ is square integrable. for almost every $u\in I$.
\end{proposition}
This proposition follows from \cite[Corollary 3.5]{tangpi-zhou22}. The well-posedness of \eqref{eqn:BSDE-graphon.decoupled} is proven in \cite[Proposition 6.1]{tangpi-zhou22} and the existence of the graphon equilibrium is established in \cite[Theorem 2.7]{tangpi-zhou22}. Thus, in view of \Cref{prop:graphon-bsde}, in order to compute the graphon equilibrium $(\opt^u)_{u\in I}$, it suffices to compute the graphon FBSDEs \eqref{eqn:BSDE-graphon.decoupled}.
Observe that unlike standard (F)BSDEs whose numerical simulations have been extensively studied (see references in the introduction), the above system of FBSDEs has two particular features: the generator depends on the law of the \emph{control} process, and
this is actually a continuum of equations all coupled through the graphon $G$. 
So far, the question of computing solutions to such equations has received little attention. In fact, even without graphon interactions, the majority of works on numerical approximations and simulations of McKean--Vlasov (F)BSDEs focus on equations with interaction through the state processes $X$ and $Y$ rather than the control process $Z$, see for example \cite{chassagneux2019numerical,angiuli2019cemracs,reisinger2020posteriori} for classical methods and \cite{fouque2020deep,carmona2022convergence,germain2022numerical,carmona2023deepfinance} for deep learning methods. 
Moreover, these works usually focus on mean-field interactions rather than graphon interactions.  To the best of the authors' knowledge, computation of graphon games has been addressed only in~\cite{aurell2022stochastic,aurell2022finite,vasal2021sequential,cuilearning}.
We are not aware of any work dealing with the question of computing solutions to continuous-space graphon games with interactions through the controls.
Notice however that \cite{aurell2022finite} considered a forward-backward system of stochastic equations with graphon interactions through the controls in the case of a finite state space. 

Before discussing our approach to the computation of the above graphon BSDE system, let us present the convergence results from the finite-agent equilibrium to the graphon equilibrium. Consider the following condition:
\begin{condition}\label{cdn:G}
 There exists a sequence of graphons $(G_N)_{N\ge1}$ such that:
   \begin{enumerate}
    \item[(1)] the graphons $G_N$ are step functions, i.e. they satisfy
       \begin{align*}
           &G_N(u,v) = G_N\Big(\frac{\ceil{Nu}}{N},\frac{\ceil{Nv}}{N}\Big) \quad \text{for}\quad (u,v) \in I \times I,\quad \text{and for every } N\in \N, 
       \end{align*}
       and it holds that
       \[N\|G_N -G\|_{2} \xrightarrow[N \to\infty]{} 0, 
       \] 
       where $\|\cdot\|_2$ is the usual $\LL^2$ norm on graphons defined as: for any graphon $G$, 
       \[ \|G\|_2 : = \bigg(\int_{I \times I} |G(u,v)|^2 \d u \,\d v \bigg)^{\frac12}.  \]
       \item[(2)]  $\lambda_{ij} = \lambda_{ji} = \mathrm{Bernoulli}(G_N(\frac{i}{N},\frac{j}{N}))$ independently for $1 \leq i,j \leq N$, and independently of $(\xi^u, \sigma^u, \theta^u, \eta^u)_{u\in I}$, and $(W^u)_{u \in I}$.
   \end{enumerate}
\end{condition}
Furthermore, let us fix some sequence of points $(u^i)_{i\ge1}$ from $I$ and re-write the strategy $\pi^i$ for agent $u^i$ in the finite-agent game to be $\pi^{i, N}$ to emphasize the dependence on the population size $N$. Then the following result holds.
\begin{theorem}
\label{thm:main.limit}
    Let \Cref{cdn: unconstrained} and \Cref{cdn:G} be satisfied. Further assume that $\E[\e^{\frac{2\rho}{\eta^i}|\xi^{i}| }]<\infty$ for all $i$ and $\xi^u \in L^2(\mu\otimes \P)$. If the $N$--agent problem \eqref{eq: n-agent obj} admits a Nash equilibrium $(\widetilde \pi^{i,N})_{i \in \{1,\dots,N\}}$, then for each $i$, up to a subsequence, the control $\widetilde\pi^{i,N}$ converges to $\widetilde \pi^{u^i}$ where $(\widetilde\pi^{u})_{u\in I}$ forms a graphon Nash equilibrium.
    In fact, we have
    \begin{align}
    \label{eq:conv.statement}
         \|\widetilde{\pi}_t^{i,N} - \widetilde{\pi}_t^{u^i}\|^2  \xrightarrow[N\to\infty]{} 0 \quad \diff t\otimes \P\text{ a.s.} \quad \text{and}\quad \Big|V^{i,N}((\widetilde{\pi}^{j,N})_{j\neq i}) - V^{\frac{i}{N},G}((\widetilde{\pi}^v)_{v\neq u^i}) \Big| \xrightarrow[N \to\infty]{} 0.
    \end{align}
\end{theorem}
This result is \cite[Theorem 2.11]{tangpi-zhou22}.
A particularly interesting case that we will often use as a benchmark to test our algorithms is when the stock price coefficients are all deterministic and strategies are unconstrained.
In this case, it turns out that for each $N\ge1$ both $N$--Nash equilibrium and graphon equilibrium exist and are given in close form. The details regarding this case are given in \cite[Proposition 2.12]{tangpi-zhou22}. For the sake of completeness, we provide the main results useful for our numerical experiments here:
\begin{proposition}\label{prop: constant coef exact solution}
    Let \Cref{cdn: unconstrained} be satisfied and let $\sigma^u$ and $\mu^u$ be deterministic measurable functions of time. Consider the following modification of the utility maximization problem \eqref{eq: n-agent obj}: $\lambda_{ii} \neq 0$, i.e, agent $i$ takes into account a weighted average of all agents' terminal wealth as their benchmark. Under this modification, we have the following results:

For all $N\in\mathbb{N}^\ast$ there is an $N$--Nash equilibrium $(\widetilde \pi^{i,n})_{i \in \{1,\dots,n\}}$ satisfying
\begin{equation*}
    \sigma^i_t\widetilde \pi^{i,N}_t = \frac{N}{N - \rho\lambda_{ii}}\eta^i_t\theta^i_t \quad \forall (N,i)\in \N^\ast\times \{1,\dots,N\}\text{ and a.e. }  t.
\end{equation*}
and a graphon Nash equilibrium $(\widetilde \pi^u)_{u \in I}$ satisfying
\begin{equation*}
    \sigma^u_t\widetilde \pi^u_t = \eta^u\theta^u_t\quad \text{a.e } (u,t)\in I\times[0,T].
\end{equation*}
Therefore, it follows that 
\begin{equation*}
    \|\sigma_t^{i}\widetilde\pi^{i, N}_t - \sigma_t^{u^i}\widetilde\pi^{u^i}_t\|  \le \frac{\rho\lambda_{ii}}{N - \rho\lambda_{ii}}\|\eta^i\theta^i\|_\infty \quad \forall (N,i)\in \N^\ast\times \{1,\dots,N\}\text{ and a.e. }  t. 
\end{equation*}
Moreover, we have the following closed-form solution for $(X^u, Y^u, Z^u)$ to the graphon FBSDEs~\eqref{eqn:BSDE-graphon}: for every $t \in [0, T]$, 
\begin{align}\label{eqn: BS-constant graphon closed form sol}
\begin{cases}
X_t^u = \int_0^t (\eta^u\theta_s^u) \cdot (\theta_s^u \d s + \d W_s^u), \\
Y_t^u = \int_t^T \bigg( \frac{\eta^u}{2} \|\theta_s^u \|^2 - \E\Big[\int_I \rho \eta^v\|\theta_s^v \|^2 G(u,v) \d v \Big]  \bigg) \d s, \\
Z_t^u = 0.
\end{cases}
\end{align}
\end{proposition}
The above results are in line with the homogeneous case studied by \citeauthor*{lacker2019mean} \cite{lacker2019mean} using PDE techniques and \citeauthor*{espinosa_touzi_2013} \cite{espinosa_touzi_2013} using BSDE techniques. When $\sigma^u$ and $\mu^u$ are constants, the stock price evolution reduces to the standard Black-Sholes model with constant coeffcients. We will present this case as the baseline case for our simulations in \Cref{sec: sanity check} and use \eqref{eqn: BS-constant graphon closed form sol} to compute the approximation errors.

\section{Description of the simulation method}

Before proceeding to the presentation of our numerical results, we first discuss the deep learning method and the neural network architecture that we employ.

\subsection{Approximate problem for graphon FBSDEs} 

We aim to compute the solutions $(X^u,Y^u,Z^u)_{u \in I}$ to the following graphon FBSDE: 
\begin{align}\label{eqn:BSDE-graphon}
\begin{cases}
    \d X_t^u = \widetilde\pi_t^u \cdot {\sigma_t^u}\{\theta_t^u \d t + \d W_t^u \},\quad 
    \opt_t^{u} = (\sigma_t^{u})^{-1}(Z_t^u + \eta^u \theta_t^u)\\
    \d Y_t^u =  \left(Z_t^u \cdot \theta_t^u + \frac{\eta^u}{2}|\theta_t^u|^2 - \E\big[\int_I \rho(Z_t^v + \eta_v\theta_t^v) \cdot \theta_t^v G(u,v) \d v \big] \right)\d t +  Z_t^u \cdot \d W_t^u  \quad \d t\otimes \d u\boxtimes\P\text{--a.s.}\\
    Y_T^u = 0, \quad X^u_0 = \xi^u.
    \end{cases}
\end{align}
In the spirit of the deep BSDE method~\cite{han2018solving} and its extension to MFGs~\cite{carmona2022convergence}, we characterize the solution to the FBSDE as the solution of an optimal control problem, and we then learn the optimal control using neural networks.  In this new optimal control problem, the goal is to find two functions, $y_0: (u,x) \in I \times \mathbb{R}^d \mapsto y_0(u,x) \in \RR$ and $z: (t,u,x) \in [0,T] \times I \times \mathbb{R}^d \mapsto z(t,u,x) \in \RR^d$ so that, by taking $Y_0^u = y_0(u,X_0^u)$ and $Z_t^u = z(t,u,X_t^u)$, we achieve $Y_T^u$ satisfying the terminal condition. 

More specifically, we consider the following graphon control problem:
\begin{equation}
\label{eq:J-y0z}
	\min_{y_0, z} J(y_0,z),
\end{equation}
where the cost function is defined by:\footnote{In our problem, the terminal cost is $0$ and hence $Y_T^u$ must be equal to $0$. It is possible to deal with more general terminal costs which depend on $X_T^u$ using the same algorithm by simply changing the cost function appearing in~\eqref{eq:J-y0z}. }
$$
	J(y_0, z) = \mathbb{E}\left[ \int_I \left| Y_T^u \right|^2 du \right]
$$
subject to: 
\begin{align}\label{eqn:BSDE-graphon-FFSDE}
\begin{cases}
    \d X_t^u = \tilde\pi_t^u \cdot {\sigma_t^u}\{\theta_t^u \d t + \d W_t^u \},\quad 
    \opt_t^{u} = (\sigma_t^{u})^{-1}(z(t,u,X_t^u) + \eta^u \theta_t^u) 
    \\
    \d Y_t^u =  \left(z(t,u, X_t^u) \cdot \theta_t^u + \frac{\eta^u}{2}|\theta_t^u|^2 - \E\big[\int_I \rho(z(t,v, X_t^v) + \eta_v\theta_t^v) \cdot \theta_t^v G(u,v) \d v \big] \right)\d t 
    \\
    \qquad \qquad  + z(t,u,X_t^u) \cdot \d W_t^u
    \\
    Y_0^u = y_0(u, X_0^u) , \quad 
    X^u_0 = \xi^u.
    \end{cases}
\end{align}

In practice, it is not possible to minimize over all functions so we replace $y_0$ and $z$ by parameterized functions, say $y_0( \cdot, \cdot; \vartheta_{y_0})$ and $z(\cdot, \cdot, \cdot; \vartheta_{z})$ with finite-dimensional parameters $\vartheta_{y_0}$ and $\vartheta_{z}$ respectively. We introduce the following cost associated to $(\vartheta_{y_0}, \vartheta_{z})$, which is the cost of using $y_0( \cdot, \cdot; \vartheta_{y_0})$ and $z(\cdot, \cdot, \cdot; \vartheta_{z})$ as control functions in~\eqref{eq:J-y0z}:
$$
	 \tilde J(\vartheta_{y_0}, \vartheta_{z}) = J(y_0(\cdot,\cdot;\vartheta_{y_0}), z(\cdot,\cdot,\cdot;\vartheta_{z})).
$$
Then the problem becomes to optimize over the parameters:
\begin{equation}
\label{minimization-graphon-over-parameters}
	\min_{\vartheta_{y_0}, \vartheta_{z}} \tilde J(\vartheta_{y_0}, \vartheta_{z}).
\end{equation}
When using deep neural networks, this is an optimization problem in finite but high dimension. To implement the computation of $\tilde{J}$, we discretize time using an Euler-Maruyama scheme, and replace the infinite population by a finite number of particles, which allows us to approximate the expectations and the integrals. 

To be specific, let us introduce a discretization of the time interval $[0,T]$ into $n^*$ sub-intervals such that $0=t_0 < t_1 < \dots < t_{{n^*}-1} < t_{n^*} = T$. For $n = 0,\dots, n^* -1$, denote the time increment by $\Delta t_n$ and the Brownian increment by $\Delta W_{t_n}$ as defined below: 
$$\Delta t_n = t_{n + 1} - t_n, \quad \Delta W_{t_n} = W_{t_{n+1}} - W_{t_n} .$$
 Furthermore, let $M$ be an integer, and consider $M$ labels $(u^i)_{i = 1, \dots M}$. These represent the (finite) population of particles used during training.
 
 Then the dynamics~\eqref{eqn:BSDE-graphon-FFSDE} is approximated by: 
\begin{align}\label{x-dynamics-discrete-time-finitepop}
X_{t_{n+1}}^{u^i} - X_{t_{n}}^{u^i} &= \opt_{t_n}^{u^i}\cdot \sigma_{t_n}^{u^i}\{\theta_{t_n}^{u^i} \Delta t_n + \Delta W_{t_n}^{u^i}  \}, \quad 
    \opt_{t_n}^{u^i} = (\sigma_{t_n}^{u^i})^{-1}(z(t_n,u^i,X_{t_n}^{u^i}; \vartheta_{z}) + \eta^{u^i} \theta_{t_n}^{u^i}).
    \\ \notag{}
    Y_{t_{n+1}}^{u^i} - Y_{t_{n}}^{u^i} &= \bigg(z(t_n,u^i, X_{t_n}^{u^i}; \vartheta_{z}) \cdot \theta_{t_n}^{u^i} + \frac{\eta^{u^i}}{2}|\theta_{t_n}^{u^i}|^2 \\ \notag
    &\qquad \qquad - \frac{1}{M}\sum\limits_{j = 1}^M\rho(z(t_n,v^j, X_{t_n}^{v^j}; \vartheta_{z}) + \eta_{v^j}\theta_{t_n}^{v^j}) \cdot \theta_{t_n}^{v^j} G(u^i,v^j)   \bigg)\Delta t_n\\ \label{y-dynamics-discrete-time-finitepop}
    &\quad + z(t_n,u^i,X_{t_n}^{u^i}; \vartheta_{z}) \cdot \Delta W_{t_n}^{u^i}, \\
    Y_0^{u^i} &= y_0(u^i, X_0^{u^i}; \vartheta_{y_0}) , \quad 
    X^{u^i}_0 = \xi^{u^i}. \label{initcond-dynamics-discrete-time-finitepop}
\end{align}
The minimization problem~\eqref{minimization-graphon-over-parameters} is replaced by:
\begin{align}\label{simulation cost function}
	\check{J}(\vartheta_{y_0}, \vartheta_{z}) = \EE \left[\frac{1}{M}\sum_{i = 1}^M|Y_T^{u^i}|^2 \right]
\end{align}

As we explain in the next subsection, we compute an approximate minimizer by using stochastic gradient descent (SGD) or one of its variants such as Adam~\cite{kingma2014adam}.

In the setting of~\eqref{simulation cost function}, one sample consists of one set of $M$ trajectories, each trajectory corresponding to one index $u \in I$. Notice that the trajectories interact through the graphon interaction term so it is not possible to sample them independently of each other. To generate $M$ trajectories, we first sample $M$ labels $u^i, i=1,\dots,M$ uniformly at random in $I$ and then simulate discrete-time trajectories using the above discrete time system~\eqref{x-dynamics-discrete-time-finitepop}--\eqref{initcond-dynamics-discrete-time-finitepop}. We use $\frac{1}{M}\sum\limits_{j = 1}^M\rho(z(t_n,v^j, X_{t_n}^{v^j}; \vartheta_{z}) + \eta_{v^j}\theta_{t_n}^{v^j}) \cdot \theta_{t_n}^{v^j} G(u^i,v^j)$ as the interaction term for $u^i$. 
Intuitively, this approximation is justified by the exact law of large numbers: for example $\E[\int_I X_t^v G(u,v) \d v] \approx \frac{1}{M} \sum_{j=1}^M X_t^{v^j} G(u, v^j)$ for $M$ large enough. We refer to~\cite{sun2006exact} for more information on the exact law of large numbers, and to~\cite{carmona2022stochastic,aurell2022stochastic,aurell2022finite} for more details on its application in the context of graphon games.

In the implementation, for the parameterized functions $y_0( \cdot, \cdot; \vartheta_{y_0})$ and $z(\cdot, \cdot, \cdot; \vartheta_{z})$, we use deep neural networks and optimize their parameters using stochastic gradient descent or one of its variants, as we explain in the next subsection.

\begin{remark}
	It is important to notice that the functions $y_0( \cdot, \cdot; \vartheta_{y_0})$ and $z(\cdot, \cdot, \cdot; \vartheta_{z})$ are not only functions of the particle's state $X_t^{u^i}$ but also of the index $u^i$. During training, we sample various indices, so that the neural networks are trained to learn the optimal control for any index. We expect this approximation to be particularly relevant when $Y^{u}_0$ and $Z^{u}_t$ depend smoothly on the index $u$. A similar approach was used in \cite{aurell2022finite} to solve graphon ODEs instead of graphon FBSDEs. We will discuss more details in \Cref{sec: losses errors and runtime} regarding the effect of the number of labels sampled $M$ on estimation accuracy.
\end{remark}

\subsection{Neural network architecture and training algorithm}

In this subsection we will discuss in details the simulation algorithm and the neural network architecture.

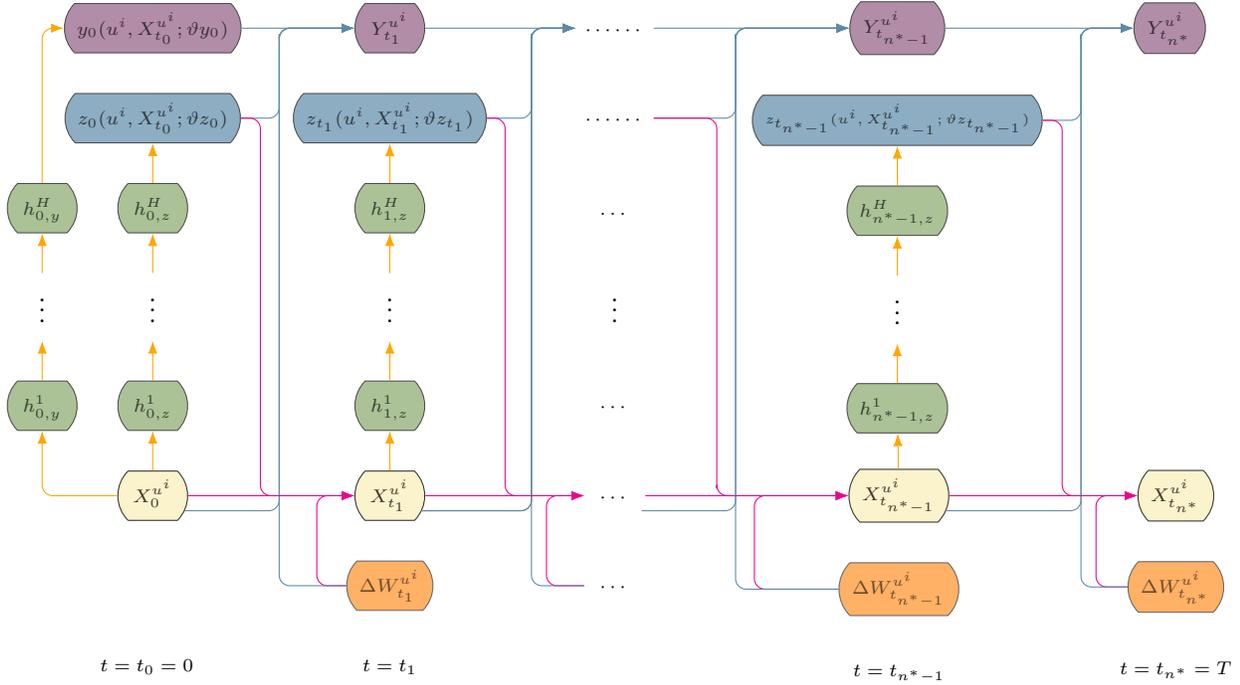
\begin{figure}[H]
\begin{tikzpicture}[scale = 0.5,
    every node/.style={font = \fontsize{7}{7}, align=center, minimum height=.8em, minimum width=.8cm,node distance=15pt}]

    \node[draw,
            rounded rectangle,
            rounded rectangle arc length = 90,
            fill = antiquefuchsia, opacity = 0.7,
            ](block00){\strut$y_0(u^{i},X_{t_0}^{u^i};\vartheta{y_0})$};
    \node[draw,
            rounded rectangle,
            rounded rectangle arc length = 90,
            fill = airforceblue, opacity = 0.7,
            below = of block00,
            ](block10){\strut$z_0(u^{i},X_{t_0}^{u^i};\vartheta{z_0})$};
    \node[draw,
            rounded rectangle,
            rounded rectangle arc length = 90,
            fill = asparagus, opacity = 0.7,
            below = of block10,
            ](block20){\strut$h_{0,z}^H$};
    \node[below = of block20]
            (block30){\strut$\vdots$};
    \node[draw,
            rounded rectangle,
            rounded rectangle arc length = 90,
            fill = asparagus, opacity = 0.7,
            below = of block30,
            ](block40){\strut$h_{0,z}^1$};

    \node[draw,
            rounded rectangle,
            rounded rectangle arc length = 90,
            fill = asparagus, opacity = 0.7,
            left = of block20,
            ](block2minus1){\strut$h_{0,y}^H$};
        \node[below = of block2minus1]
            (block3minus1){\strut$\vdots$};
    \node[draw,
            rounded rectangle,
            rounded rectangle arc length = 90,
            fill = asparagus, opacity = 0.7,
            below = of block3minus1,
            ](block4minus1){\strut$h_{0,y}^1$};

    \node[draw,
            rounded rectangle,
            rounded rectangle arc length = 90,
            fill = blond, opacity = 0.8,
            below = of block40,
            ](block50){\strut$X^{u^i}_{0}$};
    \node[draw,
            rounded rectangle,
            rounded rectangle arc length = 90,
            fill = antiquefuchsia, opacity = 0.7,
            right = 1.5cm of block00
            ](block01){\strut$Y_{t_1}^{u^{i}}$};
    \node[draw,
            rounded rectangle,
            rounded rectangle arc length = 90,
            fill = airforceblue, opacity = 0.7,
            below = of block01,
            ](block11){\strut$z_{t_1}(u^{i},X_{t_1}^{u^i};\vartheta{z_{t_1}})$};
    \node[draw,
            rounded rectangle,
            rounded rectangle arc length = 90,
            fill = asparagus, opacity = 0.7,
            below = of block11,
            ](block21){\strut$h_{1,z}^H$};
    \node[below = of block21]
            (block31){\strut$\vdots$};
    \node[draw,
            rounded rectangle,
            rounded rectangle arc length = 90,
            fill = asparagus, opacity = 0.7,
            below = of block31,
            ](block41){\strut$h_{1,z}^1$};
    \node[draw,
            rounded rectangle,
            rounded rectangle arc length = 90,
            fill = blond, opacity = 0.8,
            below = of block41,
            ](block51){\strut$X^{u^i}_{t_1}$};
   \node[draw,
            rounded rectangle,
            rounded rectangle arc length = 90,
            fill = amber(sae/ece), opacity = 0.6,
            below = of block51,
            ](block61){\strut$\Delta W^{u^i}_{t_1}$};
    \node[draw = white,
            right = 2cm of block01
            ](block02){\strut$\cdots \cdots$};
    \node[draw = white,
            below = of block02,
            ](block12){\strut$\cdots \cdots$};
    \node[draw = white,
            below = of block12,
            ](block22){\strut$\dots$};
    \node[below = of block22]
            (block32){\strut$\vdots$};
    \node[draw = white,
            below = of block32,
            ](block42){\strut$\cdots$};
    \node[draw = white,
            below = of block42,
            ](block52){\strut$\cdots$};
    \node[draw = white,
            below = of block52,
            ](block62){\strut$\cdots$};
    \node[draw,
            rounded rectangle,
            rounded rectangle arc length = 90,
            fill = antiquefuchsia, opacity = 0.7,
            right = 2.6cm of block02
            ](block03){\strut$Y_{t_{n^*-1}}^{u^{i}}$};
    \node[draw,
            rounded rectangle,
            rounded rectangle arc length = 90,
            fill = airforceblue, opacity = 0.7,
            below = of block03,
            ](block13){\tiny \strut$z_{t_{n^*-1}}(u^{i},X_{t_{n^*-1}}^{u^{i}};\vartheta{z_{t_{n^*-1}}})$};
    \node[draw,
            rounded rectangle,
            rounded rectangle arc length = 90,
            fill = asparagus, opacity = 0.7,
            below = of block13,
            ](block23){\strut$h_{n^*-1,z}^H$};
    \node[below = of block23]
            (block33){\strut$\vdots$};
    \node[draw,
            rounded rectangle,
            rounded rectangle arc length = 90,
            fill = asparagus, opacity = 0.7,
            below = of block33,
            ](block43){\strut$h_{n^*-1,z}^1$};
    \node[draw,
            rounded rectangle,
            rounded rectangle arc length = 90,
            fill = blond, opacity = 0.8,
            right = 2.7cm of block52,
            ](block53){\strut $X^{u^i}_{t_{n^*-1}}$};
   \node[draw,
            rounded rectangle,
            rounded rectangle arc length = 90,
            fill = amber(sae/ece), opacity = 0.6,
            below = of block53,
            ](block63){\strut$\Delta W^{u^i}_{t_{n^*-1}}$};
    \node[draw,
            rounded rectangle,
            rounded rectangle arc length = 90,
            fill = antiquefuchsia, opacity = 0.7,
            right = 2.5cm of block03
            ](block04){\strut$Y_{t_{n^*}}^{u^{i}}$};
    \node[draw = white,
          text = white,
            below = of block04,
            ](block14){\strut$z_{t_{n}}(u^{i},X_{t_{n}}^{u^{i}})$};
    \node[draw = white,
            text = white,
            below = of block14]
            (block24){\strut $h_{n-1,z}^1$};
    \node[draw = white,
            text = white,
            below = of block24]
            (block34){\strut $\vdots$};
    \node[draw = white,
            text = white,
            below = of block34,
            ](block44){\strut $h_{n-1,z}^1$};
    \node[draw,
            rounded rectangle,
            rounded rectangle arc length = 90,
            fill = blond, opacity = 0.8,
            right = 2.5 of block53,
            ](block54){\strut $X^{u^i}_{t_{n^*}}$};
   \node[draw,
            rounded rectangle,
            rounded rectangle arc length = 90,
            fill = amber(sae/ece), opacity = 0.6,
            below = of block54,
            ](block64){\strut$\Delta W^{u^i}_{t_{n^*}}$};

    \node[draw=white,
            below = of block61,
            ](block71){ $t = t_1$};
    \node[draw=white,
            left = 2cm of block71,
            ](block70){ $t = t_0 =  0$};
    \node[draw=white,
            below = of block63,
            ](block73){ $t = t_{n^*-1}$};
    \node[draw=white,
            below = of block64,
            ](block74){ $t = t_{n^*} = T$};

    \draw[-Latex, chromeyellow, rounded corners] (block50.west) -| (block4minus1.south);
   \draw[-Latex, chromeyellow, rounded corners] (block2minus1.north) |- (block00.west);

    \draw[-Latex, magenta, rounded corners] (block61.west) -| ($(block51.west)+(-1,-1)$) |- (block51.west);
    \draw[-Latex, airforceblue, rounded corners] (block10.east) -| ($(block51.west)+(-2,12)$) |- (block01.west);
    \draw[-Latex, magenta, rounded corners] (block10.east) -| ($(block51.west)+(-2.5,1)$) |- (block51.west);
    \draw[-Latex, airforceblue, rounded corners] (block61.west) -| ($(block51.west)+(-2,1)$) |- (block01.west);
    \draw[-Latex, airforceblue, rounded corners] ($(block50.east)+(-0.1,-.4)$) -| ($(block51.west)+(-2,1)$) |- (block01.west);

    \draw[-Latex, magenta, rounded corners] (block62.west) -| ($(block52.west)+(-1,-1)$) |- (block52.west);
    \draw[-Latex, airforceblue, rounded corners] (block11.east) -| ($(block52.west)+(-1.4,12)$) |- (block02.west);
    \draw[-Latex, magenta, rounded corners] (block11.east) -| ($(block52.west)+(-2.1,1)$) |- (block52.west);
    \draw[-Latex, airforceblue, rounded corners] (block62.west) -| ($(block52.west)+(-1.4,1)$) |- (block02.west);
    \draw[-Latex, airforceblue, rounded corners] ($(block51.east)+(-0.1,-.4)$) -| ($(block52.west)+(-1.4,1)$) |- (block02.west);    

    \draw[-Latex, magenta, rounded corners] (block63.west) -| ($(block53.west)+(-2.5,-1)$) |- (block53.west);
    \draw[-Latex, airforceblue, rounded corners] (block12.east) -| ($(block53.west)+(-3,12)$) |- (block03.west);
    \draw[-Latex, magenta, rounded corners] (block12.east) -| ($(block53.west)+(-3.5,0)$) |- (block53.west);
    \draw[-Latex, airforceblue, rounded corners] (block63.west) -| ($(block53.west)+(-3,1)$) |- (block03.west);
      \draw[-Latex, airforceblue, rounded corners] ($(block51.east)+(-0.1,-.4)$) -| ($(block52.west)+(-1.4,1)$) |- (block02.west);   
    \draw[-Latex, airforceblue, rounded corners] ($(block52.east)+(-0.1,-.4)$) -| ($(block53.west)+(-3,1)$) |- (block03.west);    

    \draw[-Latex, magenta, rounded corners] (block64.west) -| ($(block54.west)+(-1.1,-1)$) |- (block54.west);
    \draw[-Latex, airforceblue, rounded corners] (block13.east) -| ($(block54.west)+(-1.5,12)$) |- (block04.west);
    \draw[-Latex, magenta, rounded corners] (block13.east) -| ($(block54.west)+(-2,1)$) |- (block54.west);
    \draw[-Latex, airforceblue, rounded corners] (block64.west) -| ($(block54.west)+(-1.5,1)$) |- (block04.west);
    \draw[-Latex, airforceblue, rounded corners] ($(block53.east)+(-0.1,-.4)$) -| ($(block54.west)+(-1.5,1)$) |- (block04.west);

    \draw[-Latex, airforceblue]         
                  (block00) edge (block01)
                  (block01) edge (block02)
                  (block02) edge (block03)
                  (block03) edge (block04);
    \draw[-Latex, chromeyellow] 
                  (block3minus1) edge (block2minus1)
                  (block4minus1) edge (block3minus1);
    \draw[-Latex, chromeyellow] 
                  (block20) edge (block10)
                  (block30) edge (block20)
                  (block40) edge (block30)
                  (block50) edge (block40);
    \draw[-Latex, chromeyellow]
                  (block21) edge (block11)
                  (block31) edge (block21)
                  (block41) edge (block31)
                  (block51) edge (block41)
                  (block23) edge (block13)
                  (block33) edge (block23)
                  (block43) edge (block33)
                  (block53) edge (block43);
    \draw[-Latex, magenta]    
                  (block50) edge (block51)
                  (block51) edge (block52)
                  (block52) edge (block53)
                  (block53) edge (block54);
    \end{tikzpicture} 
    \caption{\small Global NN architecture for the simulation in Algorithm~\ref{alg:training-y0-z-params}. There is a sub-neural network to approximate $Y_0$ at time $0$. It corresponds to the first column. The green nodes are the intermediate layers of this network, with $h_{0,y}^\ell$ denoting the neurons in the $\ell$-th layer. This NN has $H$ layers and its parameters are denoted by $\vartheta_{y_0}$ in the text. Then, at each time step $t_n$, $n=0,\dots,n^*-1$, there is a sub-network with $H$ layers, inputs $(X_{t_n}^{u^i})_{i = 1,\dots,M}$, and outputs $z_{t_n}(u^i,X_{t_n}^{u^i}; \vartheta_{z_{t_n}}) $ as estimations for $(Z_{t_n}^{u^i})_{i = 1,\dots,M}$. There are $n^* -1$ such sub-networks, and each NN corresponds to a column above. The green nodes are the intermediate layers of one such network, with $h_{n,z}^\ell$ denoting the neurons in the $\ell$-th layer at time step $t_{n}$, for the estimation of $Z_{t_n}$. In total, there are $(H+1)(n^* -1)$ layers with free parameters to optimize. These parameters are denoted by $\vartheta_{z}$ in the text. Note that only variables that are direct outputs of neural networks are denoted using small case letters in the graph above, such as $y_0, (z_{t_n})_{n = 0,\dots, n^* - 1}$. } \label{Fig: nn-architectureFBSDE}
    \end{figure}

\paragraph{Neural network architecture.} We approximate $Y_0^{u^i}$ and $Z_0^{u^i}$ by using two $H$-layer neural networks with inputs $(u^i, X_{0}^{u^i})$, where $H$ denotes the number of layers in the neural network (NN). These two networks are characterized by the connections $X_0^{u^i} \rightarrow h_{0,z}^{1} \rightarrow \cdots \rightarrow h_{0,z}^{H} \rightarrow z_0(u^i, X_{0}^{u^i}; \vartheta_{z})$ and $X_0^{u^i} \rightarrow h_{0,y}^{1} \rightarrow \cdots \rightarrow h_{0,y}^{H} \rightarrow y_0(u^i, X_{0}^{u^i};\vartheta_{y_0})$ from \Cref{Fig: nn-architectureFBSDE}. To simplify the presentation, we use the same number of layers $H$ for all the neural networks but this number could vary from one neural network to the other. Furthermore, here and in the implementation, we used feedforward fully connected networks but other architectures could be used without changing the training algorithm described below.

\paragraph{Training method.} We aim to optimize the weights $\vartheta_{z}$, $\vartheta_{y_0}$ of these two NNs. To this end, we use the method described in Algorithm~\ref{alg:training-y0-z-params}. At iteration $k$, we compute the gradient of the loss $L(\vartheta_{y_0}^{(k)}, \vartheta_{z}^{(k)}; S^{(k)})$, (see definition in Algorithm \ref{alg:training-y0-z-params}) which depends on the sample generated at this iteration using Algorithm~\ref{alg:sampingfor-y0-z-params}, namely $S^{(k)}$, the indices $(u^i)_{i=1,\dots,M}$, initial positions $(X_0^{u^i})_{i=1,\dots,M}$ and Gaussian increments $(\Delta W_{t_n}^{u^i})_{i=1,\dots,M, n=0,\dots,n^*-1}$ sampled randomly at this iteration. Here we consider that the initial position $X_0^{u}$ of player $u \in [0,1]$ is distributed according to $\mu^{u}_0$, which may vary depending on $u$. Notice that the expectation of the empirical loss at iteration $k$ is the cost $\check{J}$ defined in~\eqref{simulation cost function}:
\[
	\EE_{S^{(k)}} [L(\vartheta_{y_0}^{(k)}, \vartheta_{z}^{(k)}; S^{(k)})] 
	= \check{J}(\vartheta_{y_0}^{(k)}, \vartheta_{z}^{(k)}),
\]
where $\EE_{S^{(k)}}$ means that the expectation is taken over the random variable $S^{(k)}$, which represents the sample picked at iteration $k$ and is of the form $S = (u^i, X_0^{u^i}, \Delta W_{t_0}^{u^i}, \dots, \Delta W_{t_{n^*-1}}^{u^i})_{i=1,\dots,M}$ with distribution $\mathcal{U}([0,1]) \otimes \mu_0^i \otimes \mathcal{N}(0, \Delta t_{0}) \otimes \cdots \otimes \mathcal{N}(0,\Delta t_{n^*-1})$. 
In practice, we do not use pure stochastic gradient descent (SGD) but rather one of its variants, namely Adam~\cite{kingma2014adam}.

\begin{algorithm}
\caption{Sampling method}\label{alg:sampingfor-y0-z-params}
\KwData{$0=t_0 < t_1 < \dots < t_{{n^*}-1} < t_{n^*} = T$: discrete time grid;  $M$: number of particles per iteration}
\KwResult{One sample $S$ consisting of labels, initial positions and Gaussian increments}
	Sample independently and uniformly at random $u^i \in [0,1]$, $i = 1, \dots M$\;
 	Sample $X_0^{u^i}$ i.i.d. according to $\mu^{i}_0$\;
	Sample Gaussian increments $\Delta W^{u^i}_{t_{n+1}} \sim \cN(0,\Delta t_n)$ for $i = 1, \dots, M$, $n=0,\dots,n^*$\; 
 	Let $S = \big(u^i, X_0^{u^i}, (\Delta W_{t_n}^{u^i})_{n=0,\dots,n^*-1}\big)_{i=1,\dots,M}$\;
	Return $S$\;
\end{algorithm}

\begin{algorithm}
\caption{Training method}\label{alg:training-y0-z-params}
\KwData{$0=t_0 < t_1 < \dots < t_{{n^*}-1} < t_{n^*} = T$: discrete time grid;  $K$: number of iterations;  $M$: number of particles per iteration; $y_0(\cdot, \cdot;\vartheta_{y_0})$ and $z(\cdot, \cdot; \vartheta_{z})$: neural networks for $Y^0$ and $Z$ with parameters denoted respectively by $\vartheta_{y_0}$ and $\vartheta_{z}$; $\vartheta_{y_0}^{(0)}$ and $\vartheta_{z}^{(0)}$: initial values for the parameters}
\KwResult{Parameters for the neural networks such that $\check{J}$ in~\eqref{simulation cost function} is (approximately) minimized}
 \For{$k=0,\dots,K-1$}{
	Using Algorithm~\ref{alg:sampingfor-y0-z-params}, generate a sample $S^{(k)} = \big(u^i, X_0^{u^i}, (\Delta W_{t_n}^{u^i})_{n=0,\dots,n^*-1}\big)_{i=1,\dots,M}$\;
	Let $Y_0^{u^i} = y_0(u^i, X_0^{u^i}; \vartheta_{y_0}^{(k)})$\;
		\For{$n = 0,\dots, n^* - 1$}{
		Let $t = t_{n+1}$\;
		Generate Gaussian increments $\Delta W^{u^i}_{t_{n+1}} \sim \cN(0,1)$ for $i = 1, \dots, M$, independently from the past and from each other\;
		Compute $Y_{t_{n+1}}^{u^i}$ and $X_{t_{n+1}}^{u^i}$ by \eqref{x-dynamics-discrete-time-finitepop} and \eqref{y-dynamics-discrete-time-finitepop} using the neural network for  current parameter $\vartheta_{z}^{(k)}$ for the neural network  $z$\;
    	} 
	Using automatic differentiation (backpropagation), compute the gradient $\nabla_{(\vartheta_{y_0}, \vartheta_{z})} L(\vartheta_{y_0}^{(k)}, \vartheta_{z}^{(k)};S^{(k)})$ of the loss $L(\vartheta_{y_0}^{(k)}, \vartheta_{z}^{(k)};S^{(k)}) = |Y_{T}|^2$\;
	Update the parameters using one step of SGD (or a variant e.g., Adam \cite{kingma2014adam}): 
	\[
		(\vartheta_{y_0}^{(k+1)}, \vartheta_{z}^{(k+1)}) = (\vartheta_{y_0}^{(k)}, \vartheta_{z}^{(k)}) -\alpha^{(k)}\nabla_{(\vartheta_{y_0}, \vartheta_{z})} L(\vartheta_{y_0}^{(k)}, \vartheta_{z}^{(k)}; S^{(k)})
	\]
 }
 Return $(\vartheta_{y_0}^{(K)}, \vartheta_{z}^{(K)})$
\end{algorithm}

\begin{remark}
We stress that one advantage of having $u$ as an input of the neural networks is that the neural networks are really functions of $u$ and hence can be used, after training, for any value of $u$, even those that have not been seen during training. 

Furthermore, here to simplify the presentation, we use one neural network per time step, but we could use a unique neural network $\vartheta_{z}$ which would take $t$ as an input, i.e., $\vartheta_{z}^{(k)}: [0,T] \times I \times \RR \to \RR$, and in the training algorithm, at iteration $k$ and time $t_n$, we would use $\vartheta_{z}^{(k)}(t_n, u^i, X_{t_n}^{u^i}; \vartheta_{z}^{(k)})$. See e.g. \cite{carmona2022convergence} for an application in the case of McKean-Vlasov dynamics. 
\end{remark}

\paragraph{Implementation detail: FBSDE vs. BSDE.} The above learning method and NN architecture were inspired by the deep BSDE method proposed by \citeauthor*{han2018solving} in \cite{han2018solving}. In terms of network structure and algorithm described above, the main difference compared to our work is the simulation of the forward processes $(X_t^{u})_{u \in I}$. In the cases that the forward process $X_t^u$ and the backward process $Y_t^u$ are truly coupled (e.g. when $\theta_t^u = X_t^u$), the model has to incorporate the Euler-Maruyama scheme for $X_t^u$ during the learning for the control $Z_t^u$ and the computation for the value process $Y_t^u$, since $X_t^u$ appears to be an input for both processes. However, in some other cases, such as when $\theta_t^u$ is a function of time $t$, or a function of the Brownian motions $W_t^u$, we notice that the forward process $X_t^u$ and the backward process $Y_t^u$ in \eqref{eqn:BSDE-graphon} decouple. In this case, the connection illustrated by the blue arrows in \Cref{Fig: nn-architectureFBSDE} for computing $Y_{t_n}^{u^i}$, simplifies from $(Y_{t_n}^{u^i}, \Delta W_{t_{n+1}}^{u^i}, X_{t_n}^{u^i}, Z_{t_n}^{u^i})\rightarrow Y_{t_{n+1}}^{u^i}$ to $(Y_{t_n}^{u^i}, \Delta W_{t_{n+1}}^{u^i}, Z_{t_n}^{u^i})\rightarrow Y_{t_{n+1}}^{u^i}$. Thus the estimation for the wealth process $X_t^u$ can be done at the end, after $Y_t^u$ and $Z_t^u$ have been learnt for each time step. This is in line with the BSDE method described in \cite{han2018solving}, and will likely reduce the computation time.

\subsection{Convergence assessment for the graphon Nash equilibrium} 
In the above problem, the goal is to make $|Y_T^u|^2$ as close to $0$ as possible. We can thus monitor the convergence of the algorithm by checking how close to $0$ the loss is. However, we also want to know whether the control learnt provides an approximate Nash equilibrium. To this end, we compute the exploitability of the learnt control. The exploitability of a given control $\pi$ is defined as the difference between the cost obtained by playing an optimal control when the rest of the population plays $\pi$ and the cost obtained by playing $\pi$ as the rest of the population. More precisely, we define the average exploitability as: 
\begin{align}\label{eqn: average explitability}
	\mathcal{E}^{N}((\pi^i)_{i=1,\dots,N}) = \frac{1}{N} \sum_{i=1}^N \mathcal{E}^{i,N}((\pi^i)_{i=1,\dots,N}),
\end{align}
where the exploitability for player $u^i$ is defined as:
$$
	\mathcal{E}^{i,N}((\pi^i)_{i=1,\dots,N}) = \sup_{\pi'} V_0^{i,G}(\pi'; (\pi^j)_{j \neq i}) - V_0^{i,G}(\pi^i; (\pi^j)_{j \neq i}).
$$
The exploitability is $0$ if and only if the set of controls is a Nash equilibrium. 
This quantity can be computed approximately using Monte Carlo samples and solving the optimization problem using a similar deep learning method. We use it to assess the convergence towards a Nash equilibrium. In the rest of this subsection, we will describe the algorithm in detail, and present the exploitability results for selected graphons and stock price models.
\paragraph{Algorithm 2.} 

\paragraph{Step 1.} Use the learning method described in \textbf{Algorithm 1} to output a trained neural network for computing solutions to the Mckean-Vlasov FBSDE \eqref{eqn:BSDE-graphon} characterizing the graphon equilibrium. Denote the trained parameters of this network by $(\vartheta_y, \vartheta_z)$.
\paragraph{Step 2.} Generate a batch of labels $(u^i)_{i = 1, \dots M}$. Obtain a batch of trajectories $(X_t^{u^i}, Y_t^{u^i}, Z_t^{u^i})_{t = t_0, t_1, \dots, t_{n^*}}$ using the trained parameters $(\vartheta_y, \vartheta_z)$.
With the batch of trajectories, calculate
\begin{itemize}
    \item equilibrium utilities $V_0^{i,G}$ for each agent $u^i$ using \eqref{eq:optim.sol.no.common} and the network outputs for $Y_0^{u^i}$.
    \item An estimate of the trajectory of the graphon mean-field $\E\big[\int_I \rho(Z_t^v + \eta_v\theta_t^v) \cdot \theta_t^v G(u,v) \d v \big]$, which appeared in \eqref{eqn:BSDE-graphon}, using $\frac{1}{M}\sum_{j = 1}^M\rho(z(t_n,v^j, X_{t_n}^{v^j}; \vartheta_{z}) + \eta_{v^j}\theta_{t_n}^{v^j}) \cdot \theta_{t_n}^{v^j} G(u^i,v^j)$ at each time step $t_n$.
\end{itemize}
\paragraph{Step 3.} Consider a new neural network with parameters $(\tilde \vartheta_{y_0}, \tilde \vartheta_z)$. Train the network to learn a new set of equilibrium strategies for players $(u^i)_{i = 1\dots M}$ using Algorithm 1, however with the graphon mean-field term in \eqref{y-dynamics-discrete-time-finitepop} replaced by the graphon mean-field calculated in \textbf{Step 2}. In other words, if we denote the output of the new network by $(\tilde X_t^{u^i}, \tilde Y_t^{u^i}, \tilde Z_t^{u^i})$, and the estimations for $(Y_0, Z)$ using the new network by $\tilde y_0, \tilde z$, the new control problem is to optimize over parameters $(\tilde \vartheta_{y_0}, \tilde \vartheta_z)$: $$
    \min_{\tilde \vartheta_{y_0}, \tilde \vartheta_{z}} \tilde J(\vartheta_{y_0}, \tilde \vartheta_{z}),
$$ subject to the new dynamics
\begin{align*}
\tilde X_{t_{n+1}}^{u^i} - \tilde X_{t_{n}}^{u^i} &= \opt_{t_n}^{u^i}\cdot \sigma_{t_n}^{u^i}\{\theta_{t_n}^{u^i} \Delta t_n + \Delta W_{t_n}^{u^i}  \}, \quad 
    \opt_{t_n}^{u^i} = (\sigma_{t_n}^{u^i})^{-1}(z(t_n,u^i,X_{t_n}^{u^i}; \boldsymbol{\tilde \vartheta_{z}}) + \eta^{u^i} \theta_{t_n}^{u^i}).
    \\ \notag{}
    \tilde Y_{t_{n+1}}^{u^i} - \tilde Y_{t_{n}}^{u^i} &= \bigg(\tilde z(t_n,u^i, X_{t_n}^{u^i}; \boldsymbol{\tilde \vartheta_{z}}) \cdot \theta_{t_n}^{u^i} + \frac{\eta^{u^i}}{2}|\theta_{t_n}^{u^i}|^2 \\ \notag
    &\qquad \qquad - \frac{1}{M}\sum\limits_{j = 1}^M\rho(z(t_n,v^j, X_{t_n}^{v^j}; \vartheta_{z}) + \eta_{v^j}\theta_{t_n}^{v^j}) \cdot \theta_{t_n}^{v^j} G(u^i,v^j)   \bigg)\Delta t_n\\
    &\quad + \tilde z(t_n,u^i,X_{t_n}^{u^i}; \boldsymbol{\tilde \vartheta_{z}}) \cdot \Delta W_{t_n}^{u^i}, \\ \notag
    \tilde Y_0^{u^i} &= \tilde y_0(u^i, X_0^{u^i}; \boldsymbol{\tilde \vartheta_{y_0}}) , \quad 
    \tilde X^{u^i}_0 = \xi^{u^i}.
\end{align*}

\paragraph{Step 4.} Compute the new equlibirum utilities $V_0^{i,G}$ obtained from the new network, and the average exploitability defined in \eqref{eqn: average explitability}

\begin{remark}
Notice that the new graphon FBSDE dynamics in \textbf{Step 3} is exactly the same as before, except that the graphon mean-field term is calculated using previously trained controls. In this case, the mean-field flow doesn't change during training, only the trajectory of players changes. Players are now playing against a frozen mean-field instead of each other. The system of FBSDE thus decouples.
\end{remark}

\section{Results}

In this section we will present various numerical results for different financial models with different graphon interactions. We will mainly consider the following Graphons
\begin{itemize}
    \item \textbf{Constant Graphon}. The graphon $G_1$ is given below:
    \begin{align*}
    G_1(u,v) = 1 \text{ for all } (u,v) \in I \times I.
    \end{align*}
    In this case, we recover a regular mean-field game in which all players are symmetric and equally interact with every other player. The mean field game is well-studied for generic games, but less well-explored for utility maximization games.
    \item \textbf{Two-block Graphon.} For some constants $a, b > 0$, The graphon $G_2$ is given by:
    \begin{align*}\
    \begin{cases}
    G_2(u,v) = a \text{ for } (u,v) \in [0,\frac12] \times [0, \frac12]\\
    G_2(u,v) = b \text{ for } (u,v) \in [\frac12, 1] \times [\frac12, 1] \\
    G_2(u,v) = 0 \quad \text{otherwise.} 
     \end{cases}
     \end{align*}
     In this case, we have two independent populations with their respective mean field games. The two populations might have different interaction strength based on the values of $a$ and $b$, but within each population players are again symmetric with the same interaction strength. A good real world example that resembles this interaction structure could be that of oil traders and electricity traders, who tend to compete amongst themselves instead of each other, despite the fact that they both trade in the energy sector.

     \item \textbf{Star Graphon}. For some constant $c > 0$ and $\alpha \in (0,1)$, The graphon $G_3$ is given by:
    \begin{align*}
    \begin{cases}
    G_3(u,v) = c \quad \text{for} \; (u,v) \in [0, \alpha] \times [\alpha, 1] \\
    G_3(u,v) = c \quad \text{for} \; (u,v) \in [\alpha, 1] \times [0, \alpha]
    \end{cases}
    \end{align*}
This is a graphon analogue of the discrete star graph. A star graph interaction in the context of the $n$-agent game models a situation where there exists one single major agent in the game who interacts with every other agent, while the other agents have no interactions among themselves. In the graphon analogue, $\alpha$-percentage of the population will be the group of major players that the rest of the population interact with. This is often the case when smaller hedge funds benchmark their performances using some average of the performances of bigger funds. Notice that in this particular case, the two groups of players (major and minor) completely disregard competitions within their own groups, and only consider the other group as benchmark.
\end{itemize}
In addition, we also consider the following two continuous-graphons
\begin{itemize}
    \item \textbf{Min-max Graphon}. $G_4(u,v) = \min(u,v)(1-\max(u,v))$ for $(u,v) \in I \times I$.
    The interaction strength becomes stronger when $u$ and $v$ are closer to $0.5$. If we consider the "closeness" between label $u$ and $v$ as a measure of similarity between players, the min-max graphon models the following interaction structure: when considering other funds as performance benchmarks, a fund puts more weight onto another fund that is more "similar" to itself and less weight onto the performance of a less "similar" fund. Notice that here similarity can carry a range of different definitions, such as fund size, fund type, major products they trade in, etc.
    \item \textbf{Power-law Graphon}. $G_5(u,v) = (uv)^{-\gamma}, \gamma \in \R$. This graphon has applications in coupled dynamical systems, see e.g. \cite{medvedev2017kuramoto}. Readers can refer to \cite{aurell2022stochastic} for detailed discussion regarding the power-law graphon and its application in linear quadratic stochastic graphon games. For this work, we will fix $\gamma = -\frac12$. Given $\gamma = -\frac12$, the interaction strength will be stronger when $u,v$ are closer to $1$.
\end{itemize}
\subsection{Sanity check: Black-Scholes model with constant coefficients.}\label{sec: sanity check}
We first present simulation results under the standard Black-Scholes market model where $\sigma_t^u$ and $\theta_t^u$ are constants in \eqref{stock evolution}. Notice that our assumption does not reduce the model to the case where all players invest in a single stock, since each stock $S^u_t$ has their own respective Brownian motion $W_t^u$ driving the price. In addition to constant coefficients for the stocks, we further make the assumption that all players have the same risk preference parameter. In other words, we assume that 
\begin{align}\label{assumption bs}
    \sigma_t^u = \sigma, \;\theta_t^u =\theta,\; \text{and } \eta^u = \eta \text{ for all } u \in [0,1] \text{ and } t \in [0, T].  
\end{align}
\subsubsection{Path $Y_t^u$ and utility under graphon equilibrium for different graphons}
\begin{figure}[H]
\minipage{0.33\textwidth}
    \includegraphics[width = \linewidth]{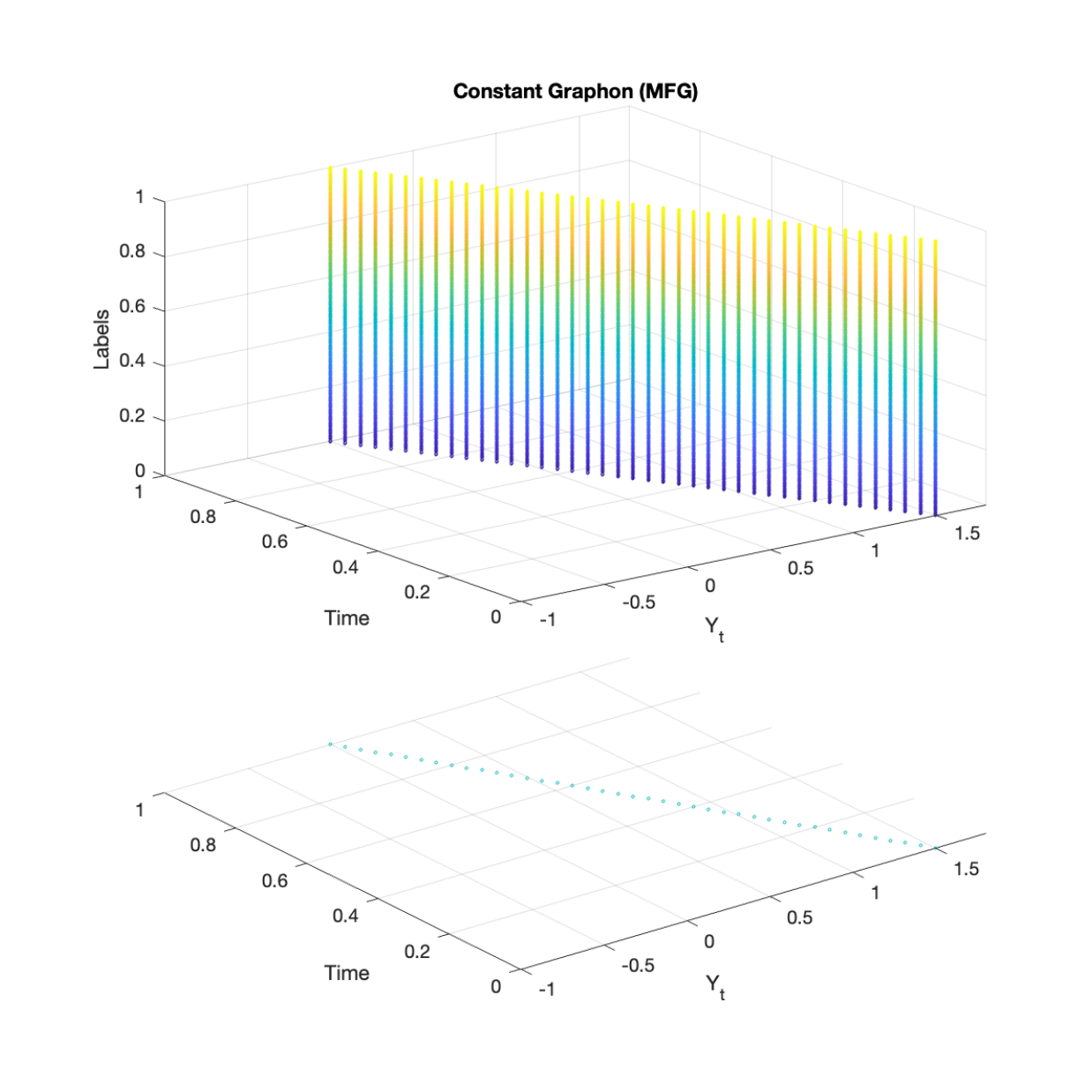}
\endminipage
\minipage{0.33\textwidth}
    \includegraphics[width = \linewidth]{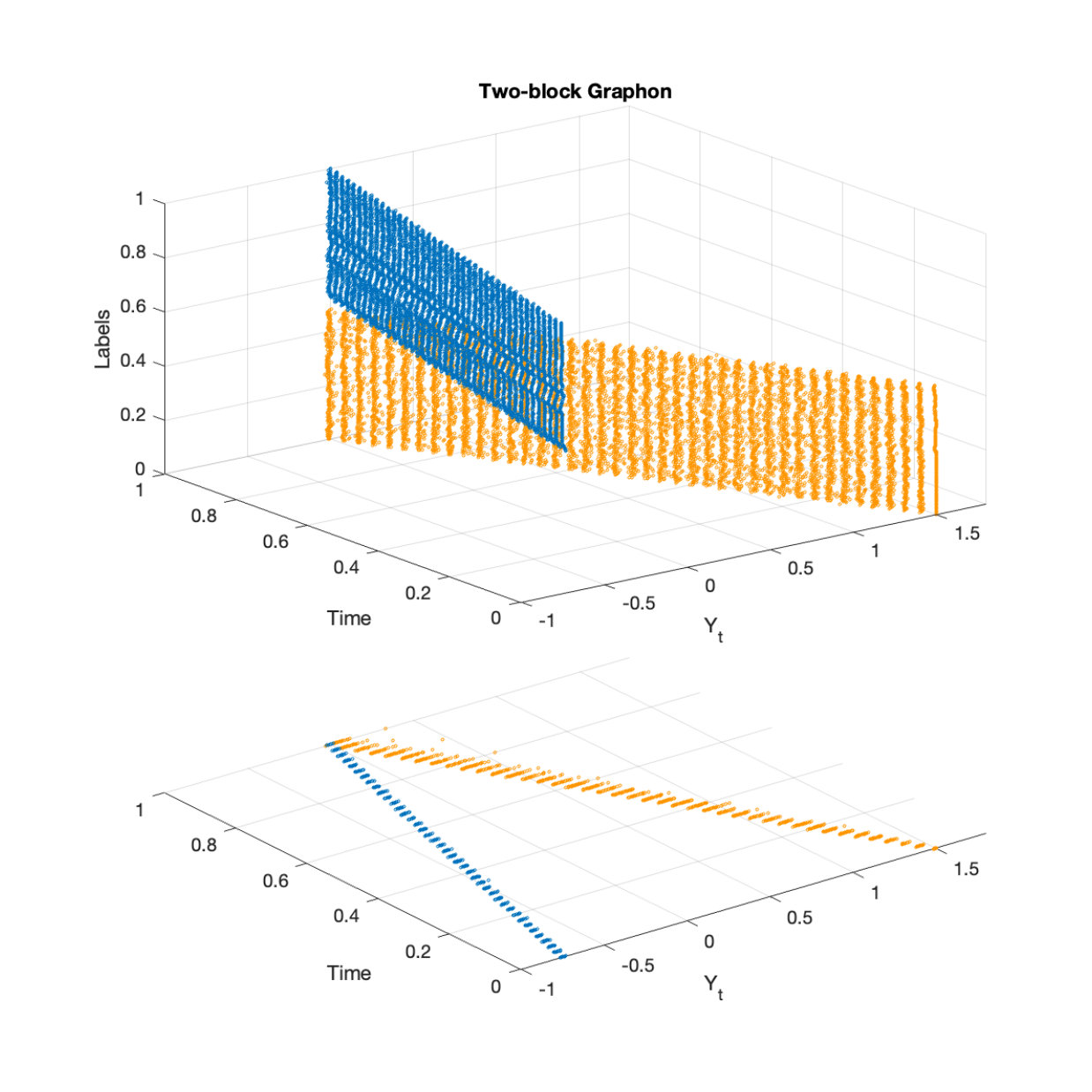}
\endminipage
\minipage{0.33\textwidth}
    \includegraphics[width = \linewidth]{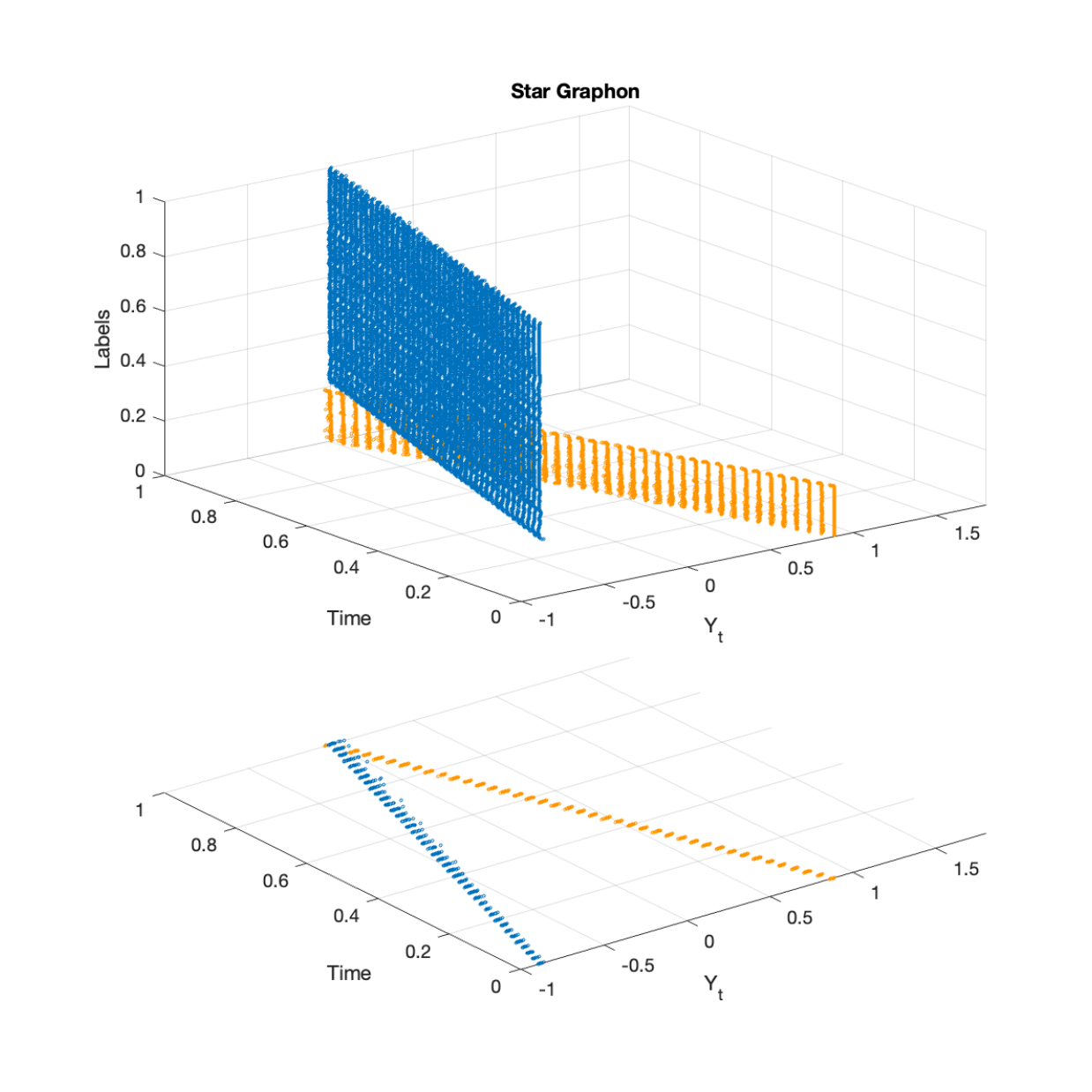}
\endminipage

\caption{\textbf{Top:} The trajectory of the value function $Y_t$ against time $t$ for different labels $u \in I$. \textbf{Bottom:} The projection of the top panels onto the $(t, Y_t)$ plane. For the constant graphon case, darker colors correspond to indices closer to $0$ and lighter colors correspond to indices closer to $1$. For the two-block graphon, the orange trajectories corresponds to the population with indices in $[0, 0.5)$, and the blue trajectories correspond to the population with indices in $[0.5, 1]$. We chose $a = 2$ and $b = 0.5$ for $G_2$. For the star graphon, we chose $\alpha = 0.2$. Thus the orange trajectories corresponds to the $20\%$ of the population that are major players, and the blue trajectories corresponds to the $80\%$ of the population that are minor players. }\label{Figure: constant BS piecewise}
\end{figure}

\begin{figure}[H]
\minipage{0.49\textwidth}
    \includegraphics[width = 0.8 \linewidth, right]{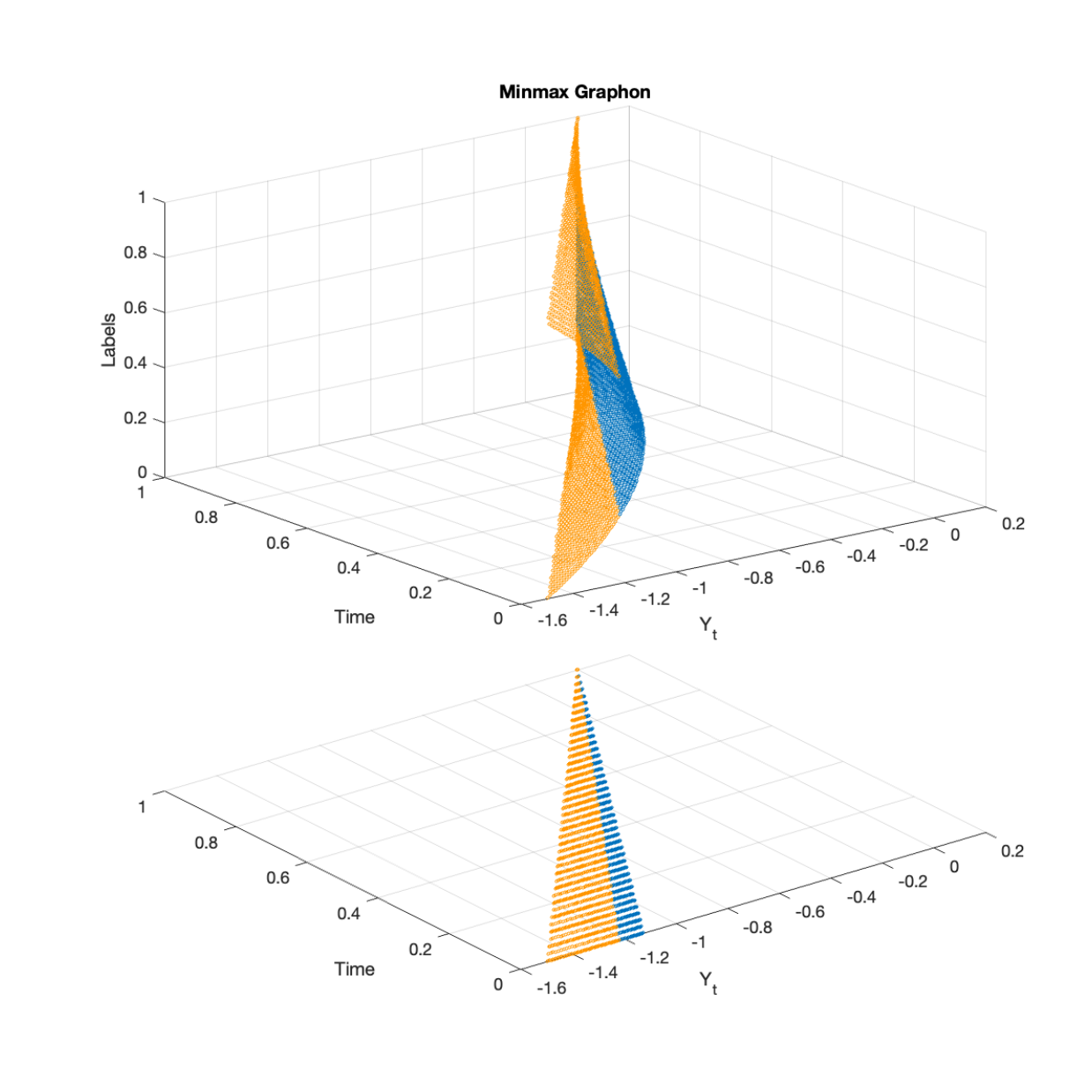}
\endminipage
\minipage{0.49\textwidth}
    \includegraphics[width = 0.8 \linewidth, left]{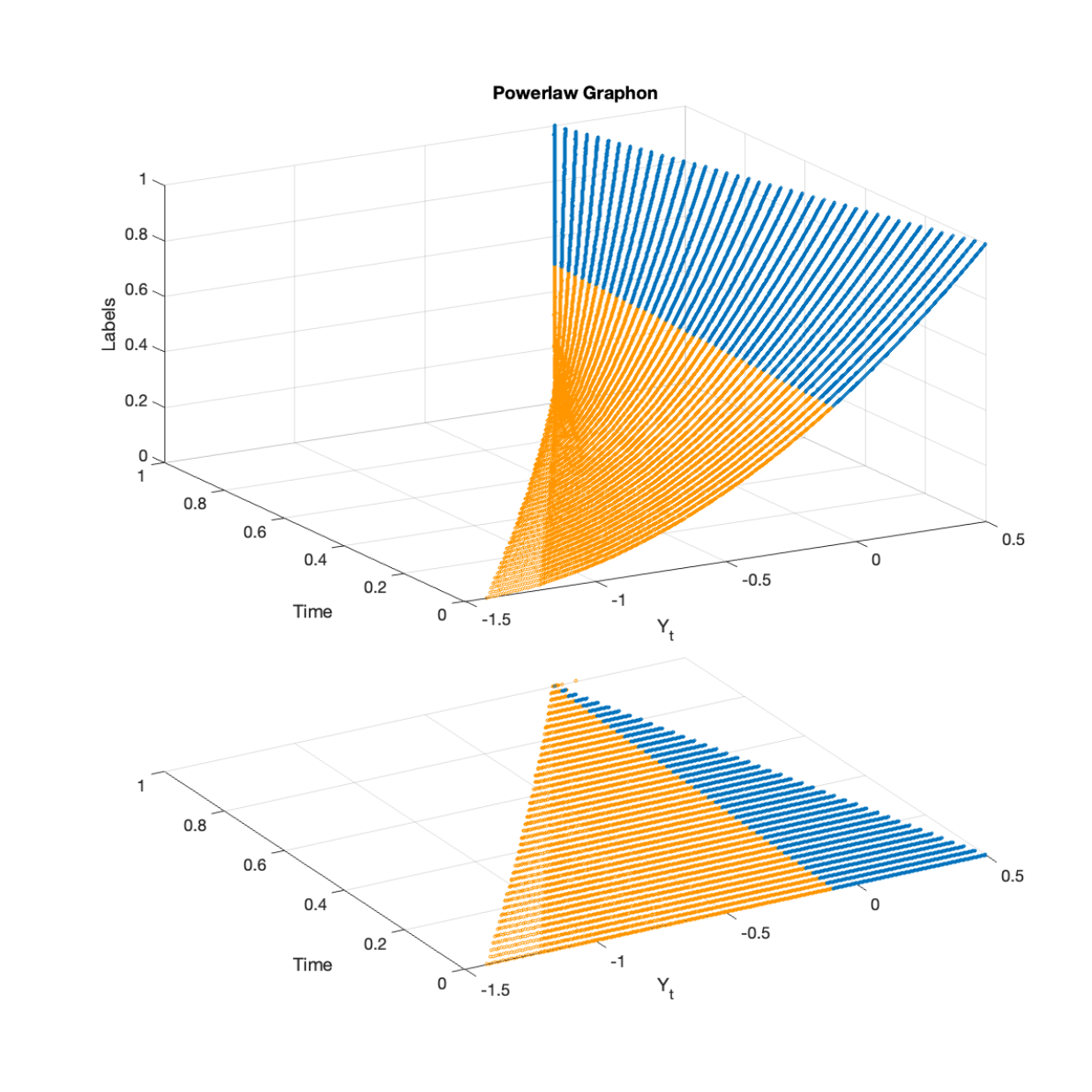}
\endminipage

\caption{\textbf{Top:} The trajectory of the value function $Y_t$ against time $t$ for different labels $u \in I$. \textbf{Bottom:} The projection of the top panels onto the $(t, Y_t)$ plane. For the min-max graphon, the blue trajectories correspond to the population with indices closer to 0.5. The orange trajectories correspond to the population with indices further away from 0.5 For the power-law graphon, the orange trajectories corresponds to the population with indices in $[0, 0.5)$, and the blue trajectories correspond to the population with indices in $[0.5, 1]$.}\label{Figure: constant BS continuous}
\end{figure}

In light of \Cref{prop: constant coef exact solution}, we have exact solutions for $(Y_t, Z_t)$ in this case. In particular, $Z_t = 0$ for all $t \in [0,T]$, and $Y_t$ reduces to linear functions of time. 
The exact solutions for $Y_t$ corresponding to the three piecewise-constant graphons are: 
\begin{itemize}
    \item When $G = G_1$, $Y_t^u = (\rho - \frac12) \eta(T - t)$ for all $u \in I$.
    \item When $G = G_2$, $Y_t^u = \frac{\eta}{2}(\rho a - 1)(T-t)$ for $u \in [0,\frac12)$ and $Y_t^u = \frac{\eta}{2}(\rho b - 1)(T-t)$ for $u \in [\frac12, 1]$.
    \item When $G = G_3$, $Y_t^u = (1-\alpha)\rho\eta - \frac{\eta}{2}(T-t)$ for $u \in [0, \alpha)$ and $Y_t^u = \alpha \rho \eta - \frac{\eta}{2}(T-t)$ for $u \in [\alpha, 1]$.
\end{itemize}
We simulated the graphon equilibriums for all five graphons mentioned above with coefficient values $\sigma = 0.1, \theta = 1, \eta = 3, \rho = 1$. We also chose the terminal time $T = 1$, and divided the time interval into $N = 40$ equidistant subintervals. The simulation results for the path $(Y_t^u)_{t \in [0,T]}$, $u \in I$ are shown in \Cref{Figure: constant BS piecewise} below for piecewise constant graphons $G_1$, $G_2$ and $G_3$, and in \Cref{Figure: constant BS continuous} below for continuous graphons $G_4$ and $G_5$. 

From the paths $Y_t^u$ given in \Cref{Figure: constant BS piecewise} above, we can see that for a mean-field game, $Y_t^u$ for different $u \in I$ are indistinguishable since all players are indistinguishable and $Y_t^u$ is a deterministic process. Both the two block graphon and the star graphon separate the players into two groups, resulted in two major paths on the projection plane. For the two-block graphon, the two groups of players are independent and the trajectory of $Y_t^u$ depends on the specific interaction strength of the group that player $u$ belongs to. For the star graphon, the two groups of players do affect each other, and the trajectory is solely dependent upon the scale of the other group.
A higher interaction strength for $G_2$, or a bigger population size for the interacting group (i.e., the other group) for $G_3$ both leads to higher values of $Y_0^u$ and lower values of the equilibrium utility $V_0^{u,G}$. These phenomenons are more clearly demonstrated in \Cref{Figure: utility vs. labels} below.

From the paths $Y_t^u$ given in \Cref{Figure: constant BS continuous} below, it is clear to see that for a min-max graphon, players with indices closer to $0.5$ obtain higher values of $Y_0$ and lower equilibrium utilities, similarly for players with indices closer to $1$ for the power-law graphon case.

In light of \Cref{prop:graphon-bsde}, we can calculate the graphon equilibrium utilities for any $u \in I$ given the $Y_0^u$ obtained from the network. The figure below plots $V_0^{u,G}$ as a function of $u \in I$ for five different graphons. Note that all our simulations were run using the same value for the starting wealth $X_0^u$ for every $u \in I$. Given the same amount of starting wealth, risk aversion parameter $\eta$ and price coefficient $\theta$,  one can see from \eqref{eq:optim.sol.no.common} and \eqref{eqn: BS-constant graphon closed form sol} that the equilibrium utility is entirely dependent upon the interaction strength $G(u,v)$, as evidenced in \Cref{Figure: utility vs. labels} below. Higher interaction strength leads to lower equilibrium utility, which can be explained by the fact that more interactions mean more competition. In the case of $G_2$, given similar assets and risk-aversion attitudes, more competition leads to lower utilities. However in the case of $G_3$, this result might seem counter-intuitive, since the "major" group of players ended up with lower utilities. Keep in mind that in a real world scenario, the major players usually own a larger amount of capital and are likely to start with higher initial wealth $X_0^u$, which will lead to higher values of terminal utility as reflected in \Cref{eq:optim.sol.no.common}.
\begin{figure}[H]
\centering
\includegraphics[width = .5\textwidth]{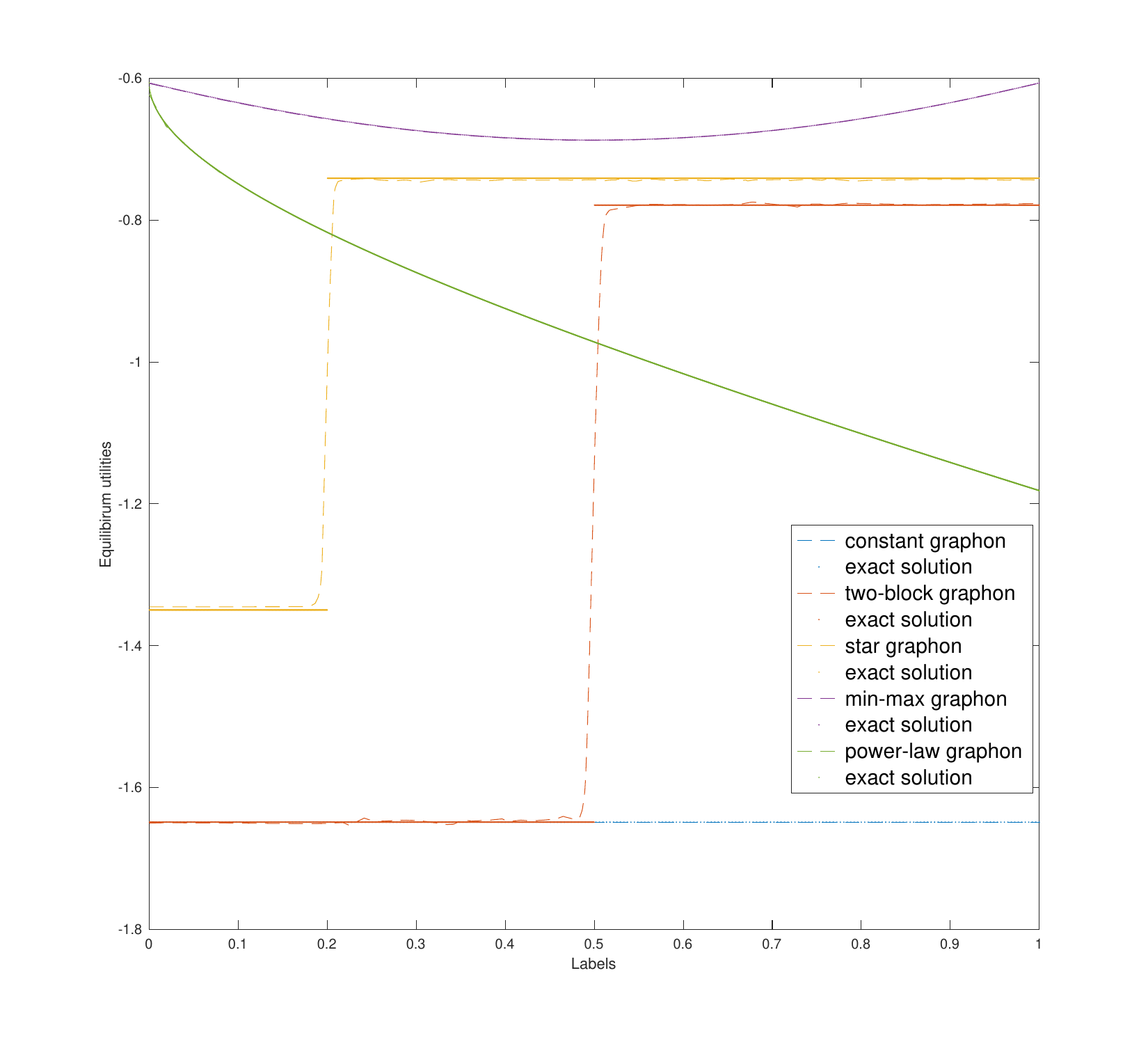}
\caption{Equilibrium utilities vs. labels for different graphons}
\label{Figure: utility vs. labels}
\end{figure}

\subsubsection{Losses, errors and runtime}\label{sec: losses errors and runtime}
Given the existence of exact solutions, we can also compute the relative approximation error for $Y_0$. For $G_1$, we achieved a validation loss on the order of $10^{-11}$ and a relative approximation error on the order of $10^{-7}$ percent in a runtime of $408$ seconds. The results are shown in the two panels on the right below. The shaded band in each graph indicates a $95\%$ confidence interval with $30$ different runs.
\begin{figure}[H]
\minipage{0.49\textwidth}
    \includegraphics[width = 0.9 \linewidth, right]{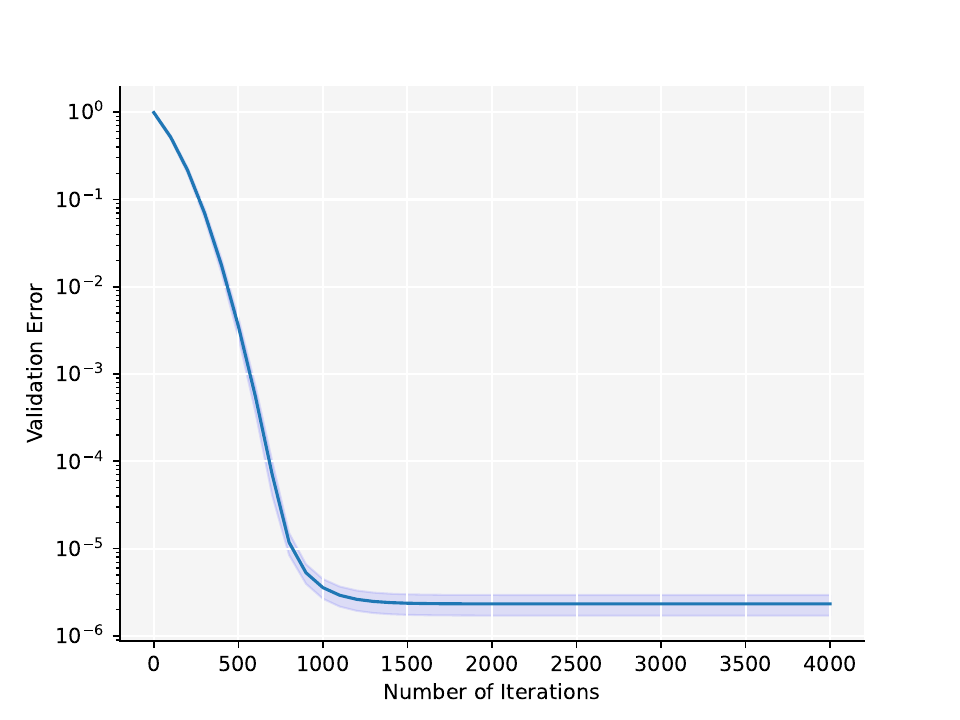}
\endminipage
\minipage{0.49\textwidth}
    \includegraphics[width = 0.9 \linewidth, left]{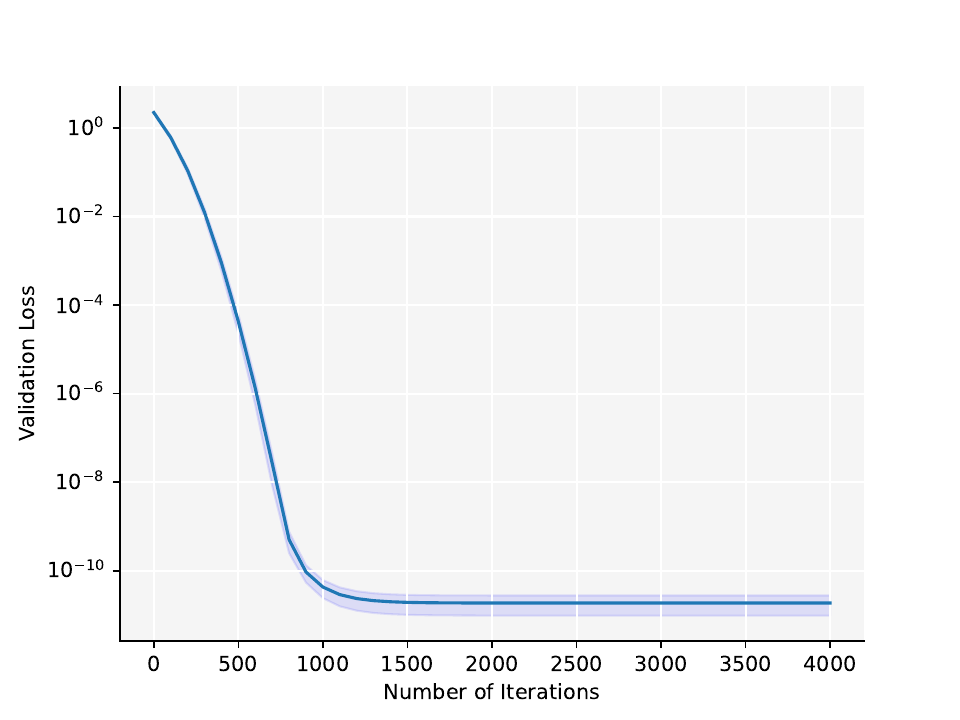}
\endminipage

\caption{Validation errors and validation losses for $G_1$ }\label{Figure: error and validation loss}
\end{figure}

\begin{figure}[H]
\minipage{0.49\textwidth}
    \includegraphics[width = 0.95 \linewidth, right]{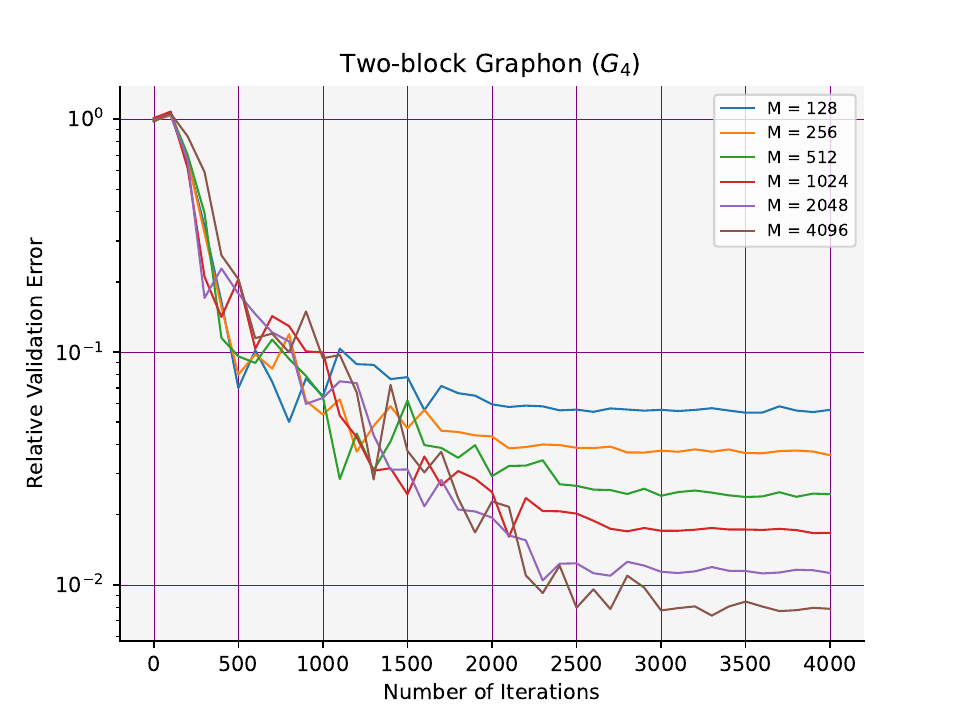}
\endminipage
\minipage{0.49\textwidth}
    \includegraphics[width = 0.95 \linewidth, left]{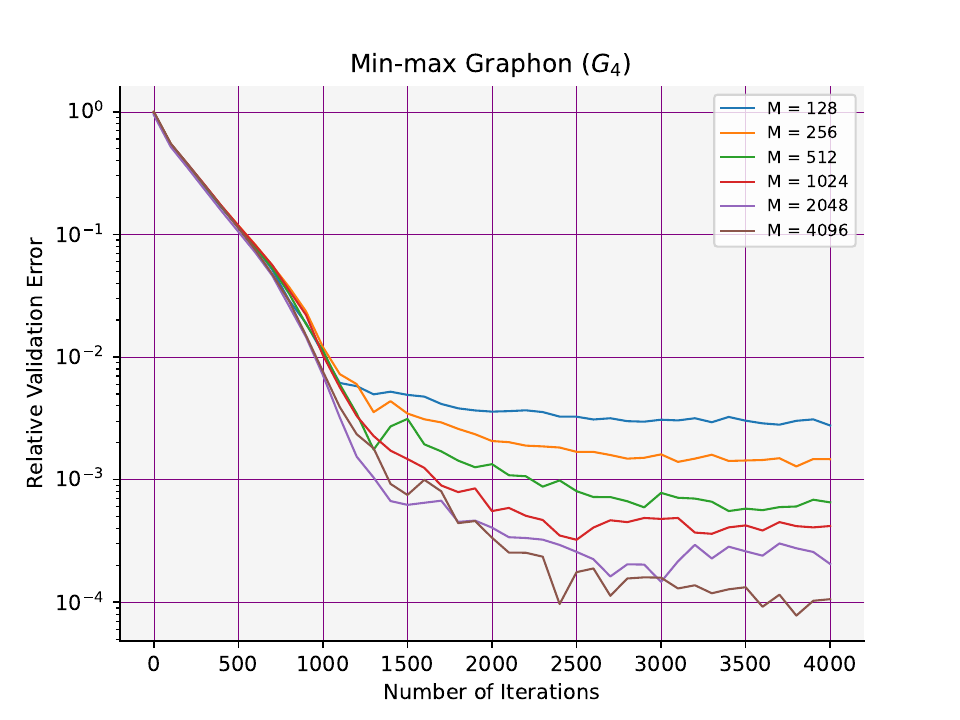}
\endminipage

\caption{{}Validation errors for different values of $M$ for $G_2$ and $G_4$}
\label{Figure: error vs. dimension}
\end{figure}

While simulating the above graphon games, for non-constant graphons, we notice a trade-off between approximation errors and runtime when we increase the number of trajectories $M$ per iteration. Intuitively, $M$ indicates the size of the population we're using to approximate the graphon game. Meanwhile, a higher $M$ also gives the network more information about the graphon in the sense of providing more data points on a continuous (or piece-wise constant) function. We ran simulations for different values of $M$ and using the two-block graphon $G_2$ and the min-max graphon $G_4$. The relative validation errors in \Cref{Figure: error vs. dimension} and the runtime in \Cref{Table: runtime vs. dimension} are computed as averages over 8 different runs. As illustrated in \Cref{Figure: error vs. dimension} and \Cref{Table: runtime vs. dimension}, higher $M$ leads to lower approximation error at the cost of longer runtime.

\begin{table}[H]
\centering
    \begin{tabular}{c|c c c c c c}
\hline
number of trajectories per batch (M) & 128 & 256 & 512 & 1024 & 2048 & 4096 \\
\hline\hline
runtime for $G_2$ (minutes) & 5.08 & 6.38 & 8.71 & 13.78 & 23.94 & 44.85 \\
\hline
runtime for $G_4$ (minutes) & 5.19 & 6.44 & 8.83 & 13.87 & 24.05 & 44.95 \\
\hline
\end{tabular}
    \caption{Runtime for different values of $M$ for $G_2$ and $G_4$}
    \label{Table: runtime vs. dimension}
\end{table}

\subsection{Black-Scholes model with Markovian coefficients.}
\label{sec:BS.Markov} For this specific case, we assume that
\begin{align}\label{assumption bs markovian}
    \sigma_t^u = \sigma, \;\theta_t^u =W_t^u,\; \text{and } \eta^u = \eta \text{ for all } u \in [0,1] \text{ and } t \in [0, T].  
\end{align}
 The simulation results for the path $(Y_t^u)_{t \in [0,T]}$ are shown in \Cref{Figure: markovian BS} in the appendix. In the case that $\theta_t^u = W_t^u$, we can see that the path $Y_t^u$ are stochastic. The separation of the paths into different groups based on the values of $u$ are similar as illustrated in \Cref{Figure: constant BS piecewise} and \Cref{Figure: constant BS continuous} for the constant coefficient case under the assumption \eqref{assumption bs}.

We are interested in simulating and learning the expected equilibrium wealth and equilibrium utility for different players $u \in I$ with different graphon interactions, since in this case there are no known analytic solutions to these quantities. we consider simulating the path of the expected wealth, and expected benchmark wealth for a particular agent $u \in I$ in a graphon equilibrium, i.e., the following two quantities:
\begin{align*}
    \E[X_t^u] \quad \text{and} \quad \E[X_t^u - \int X_t^v G(u,v) \d v], \quad \text{for} \;t \in [0,T].
\end{align*}
We conducted experiments by taking the interaction strength $c = 1$ and the proportion of major players $\alpha$ to be $0.5$, $0.2$, $0.1$ and $0.05$ respectively. Notice that in the mean-field graphon and the star graphon cases, the players are indistinguishable within their respective "group". Hence we can approximate the above two quantities by taking the average of wealth across labels within each group. This average is then representative of the expected wealth of a single player in that group. We plotted the two quantities above for the star graphon, as well as for the mean-field graphon $G_1$ as a benchmark for comparison. 

\begin{figure}[H]
\minipage{0.49\textwidth}
    \includegraphics[width = \linewidth, right]{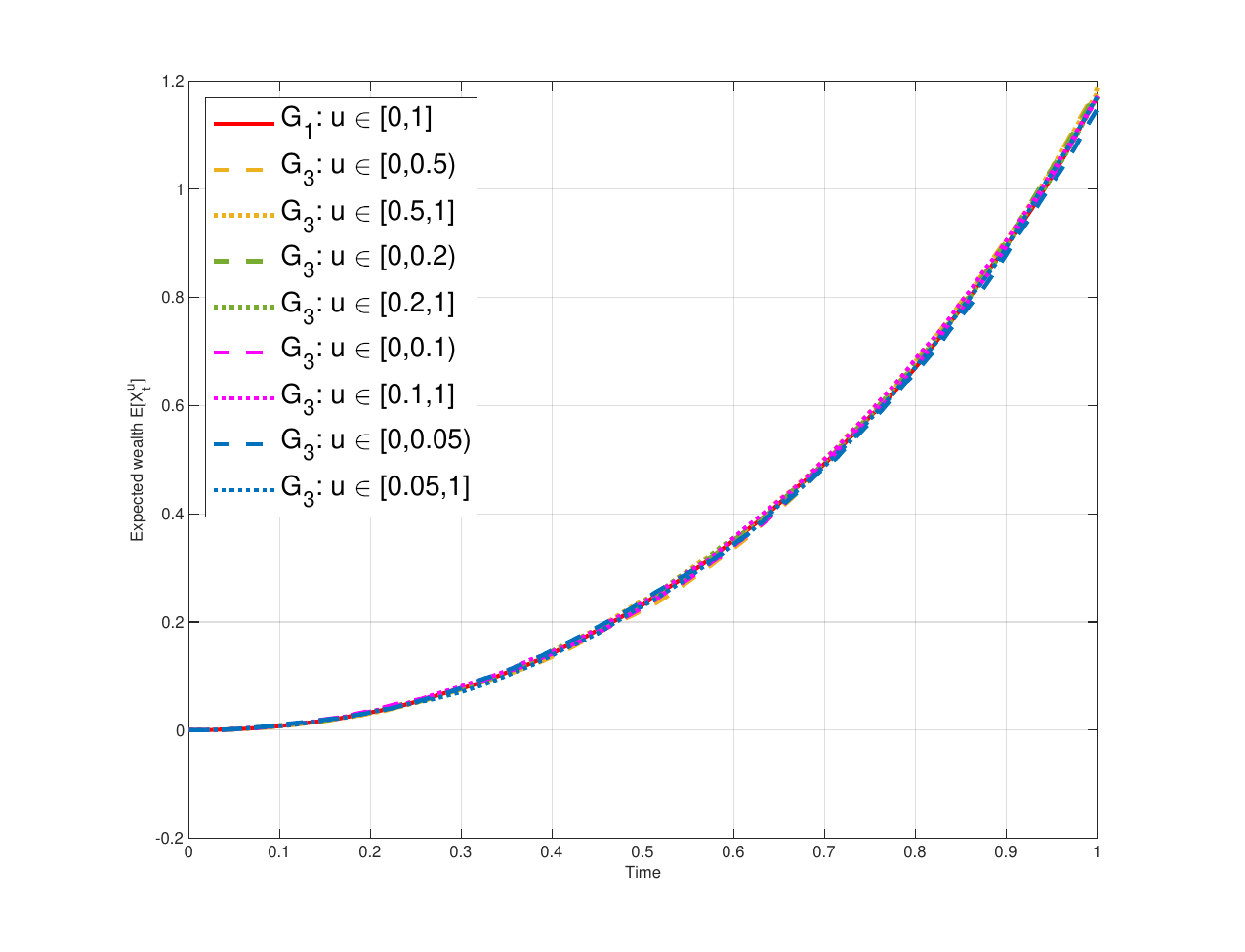}
\endminipage
\minipage{0.49\textwidth}
    \includegraphics[width = \linewidth, left]{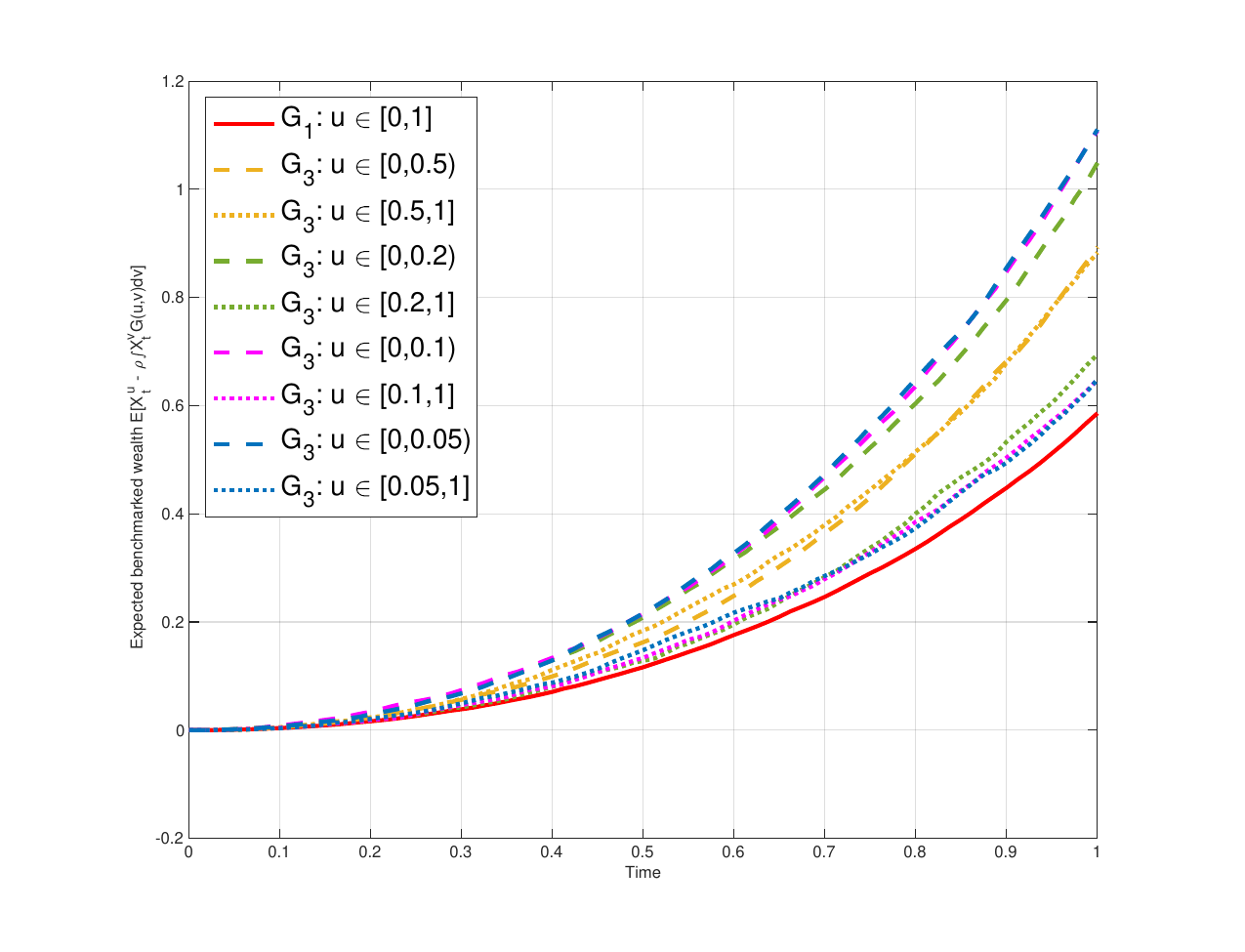}
\endminipage

\caption{Expected wealth and benchmarked wealth over time when $\eta^u = \eta$}
\label{Figure: expected wealth and relative wealth - same eta}
\end{figure}

The figure on the right indicated that for the star graphon, expected benchmarked wealth for the major player group increases, and the same quantity for the minor player group decreases, as proportion of major player $\alpha$ decreases. The expected benchmark wealth over time in the case of a star graphon are in general higher than the same quantity for players in the case of a mean-field game, with this quantity for the minor players in the star graphon case decreasing and approaching the same quantity in the case of a mean-field game as $(1-\alpha)$ increases. It is worth noticing from the figure in the left that although the expected benchmark wealth differs from players in different "groups", the expected wealth for the two groups coincides regardless of the value of $\alpha$, and coincides with the same quantity for players in the mean-field game. In other words, $\E[X_t^u]$ is independent of the label $u$. This result follows from the following proposition:
\begin{proposition}
\label{prop: wealth independent of label}
Consider the following Mckean-Vlasov type (F)BSDEs characterizing a graphon equilibrium as described in \Cref{prop:graphon-bsde}, with $\sigma_t^u = \sigma >0$, $\eta_t^u = \eta > 0$, and $\theta_t^u = W_t^u$ for all $(t,u) \in [0,T] \times I$:
\begin{align}\label{eqn:BSDE-graphon-markovian}
\begin{cases}
    \d X_t^u = \tilde\pi_t^u \cdot {\sigma}\{W_t^u \d t + \d W_t^u \},\quad 
\opt_t^{u} = (\sigma)^{-1}(Z_t^u + \eta W_t^u)\\
    \d Y_t^u =  \left(Z_t^u \cdot W_t^u + \frac{\eta}{2}|W_t^u|^2 - \E\big[\int_I \rho(Z_t^v + \eta W_t^v) \cdot W_t^v G(u,v) \d v \big] \right)\d t +  Z_t^u \cdot \d W_t^u  \quad \mu\boxtimes\P\text{--a.s.}\\
    Y_T^u = 0, \quad X^u_0 = \xi^u.
    \end{cases}
\end{align} 
If a solution $(X^u, Y^u, Z^u)$ exists for \eqref{eqn:BSDE-graphon-markovian} for $t \in [0.T]$. Then it holds that $\E[Z_t^u]$ is independent of $u$ for $t \in [0,T]$.
\begin{proof}  
The easiest way to prove the result uses tidbits of the theory of Malliavin calculus.
We refer the reader for instance to \cite{delong2010malliavin} for elements of the theory applied to BSDEs.
Below, we use the notation $(D^u_t\xi)_{t\in [0,T]}$ for the Malliavin derivative of a random variable $\xi$ in the direction of the Brownian motion $W^u$.
From standard results for FBSDEs  (see e.g. \cite{delong2010malliavin}), the solution $(Y_t^u, Z_t^u)_{t \in [0,T]}$ is differentiable in Malliavin's sense, and there exists a version of the Malliavin derivative (in the direction of Brownian motion $W^u$) denoted by $(D_t^uY_t^u)_{t \in [0,T]}$ which satisfies $D_t^uY_t^u = Z_t^u$.
Moreover, the Malliavin derivatives $(D^u_\tau Y^u, D^u_\tau Z^u)$ satisfy the dynamics
\begin{align*}
    D^u_\tau Y^u_t&= \int_t^T \left( D_\tau^uZ_s^u \cdot W_s^u + Z_s^u + \eta W_s^u \right)\d s - \int_t^T D_\tau^uZ_s^u \cdot \d W_s^u,\quad \tau \le t \\
\quad &= \int_t^T \big( Z_s^u + \eta W_s^u \big) \d s - \int_t^T    D_\tau^uZ_s^u \cdot (\d W_s^u - W_s^u \d s)
\end{align*}
and $(D^u_\tau Y^u_t, D^u_\tau Z^u_t) = (0,0)$ for $\tau>t$.
Let us introduce $\Q$, the probability measure with density $\frac{\d \Q}{\d \P} := \cE (\int_0^T W_s^u \d W^u_s) := \exp(\int_0^TW_s^u\d W^u_s - \frac12\int_0^T|W^u_s|^2\d s)$.
The process $W_t^u - \int_0^t W_s^u \d s$ is a standard Brownian motion under $\Q$. 
Recalling that $D^u_tY^u_t = Z^u_t$ and taking conditional expectation under $\Q$, we have that 
    $$Z_t^u = \E^{\Q}\bigg[ \int_t^T \big( Z_s^u + \eta W_s^u \big) \d s \big| \cF_t \bigg].$$
and thus that $Z^u_t = e^{T-t}\eta\E^{\Q}[\int_t^TW_s^u\d s\big| \cF_t]$.
Hence, taking the expectation on both sides leads to
\begin{align*} 
    \E[Z_t^u] &= e^{T-t}\eta \E\bigg[\cE\bigg(\int_t^TW^u_s\d W^u_s\bigg)\int_t^TW^u_s\d s \bigg]\\
              =&e^{T-t}\eta \E\bigg[\cE\bigg(\int_t^TW^v_s\d W^v_s\bigg)\int_t^TW^v_s\d s \bigg]
\end{align*}
where the second equality follows from the fact that $W^u$ and $W^v$ have the same distribution.
This concludes the proof.
\end{proof}
\end{proposition}
Here since we take the risk-aversion parameter $\eta^u$ to be constant across all players $u \in I$, by \Cref{prop: wealth independent of label}, the expected wealth over time for players with labels from different groups coincide. One might naturally wonder how the expected wealth and terminal utilities under equilibrium will behave when $\eta^u$ is non-constant. These discussions are presented in the following section.

\subsubsection{Effect of the risk-aversion parameter $\eta^u$ on equilibrium wealth and utilities}

We explore the effect of $\eta^u$ on $\E[X_t^u]$ and $\E[X_t^u - \int X_t^v G(u,v) \d v]$ for the star graphon and the minmax graphon. For the star graphon, we used $\alpha = 0.5$ and $\eta^u = 0.75$ as a benchmark, such that both quantities mentioned above coincide for the two groups, as illustrated by the red lines in the figures below. We experimented with choosing $\eta^u$ to be a constant higher than $0.75$ for $u \in [0, 0.5)$, and a constant lower than $0.75$ for $u \in [0.5. 1]$. The figures indicated that smaller risk aversion parameters lead to higher expected wealth and higher expected benchmarked wealth over time. 

\begin{figure}[H]
\minipage{0.49\textwidth}
    \includegraphics[width = \linewidth, right]{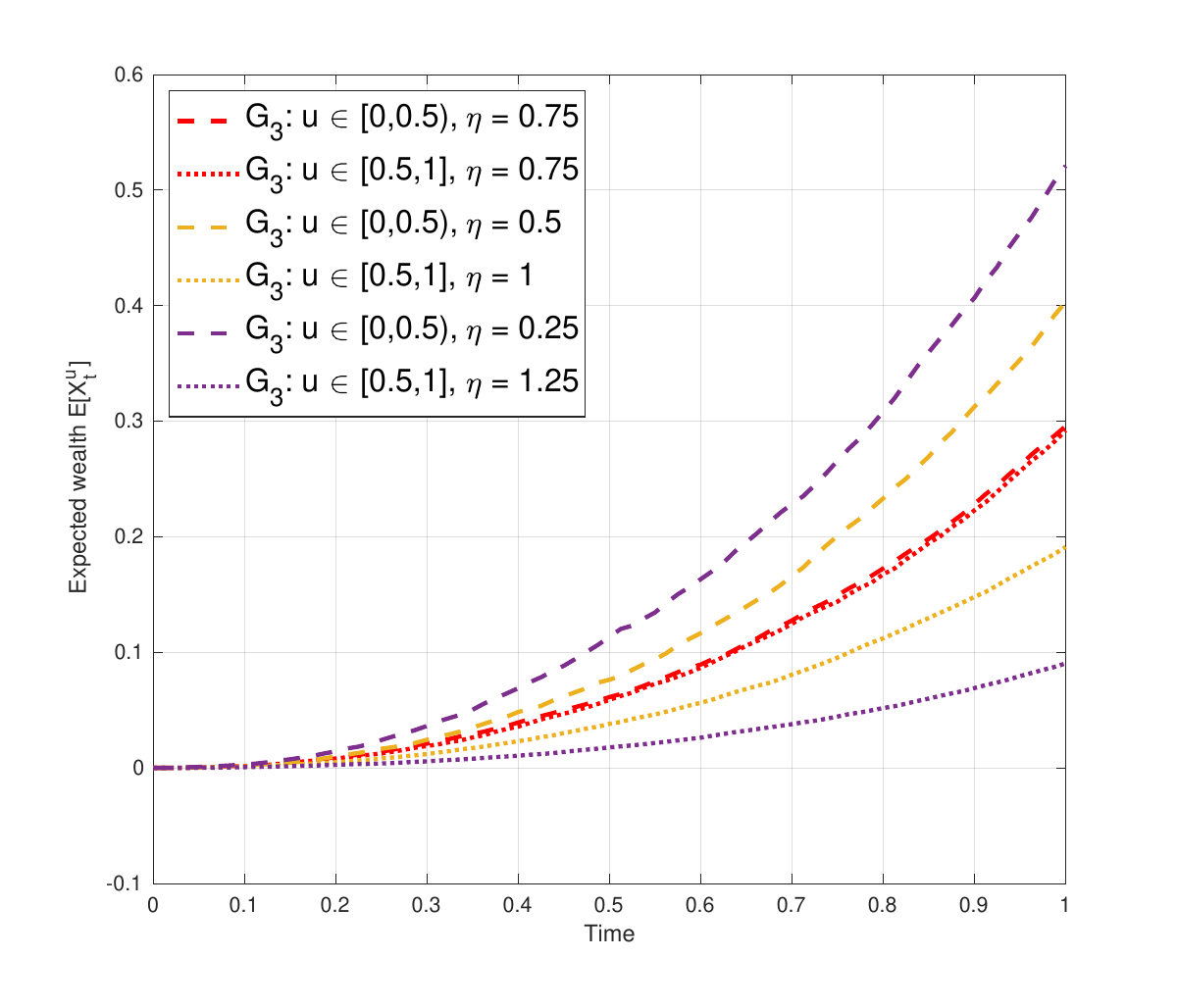}
\endminipage
\minipage{0.49\textwidth}
    \includegraphics[width = \linewidth, left]{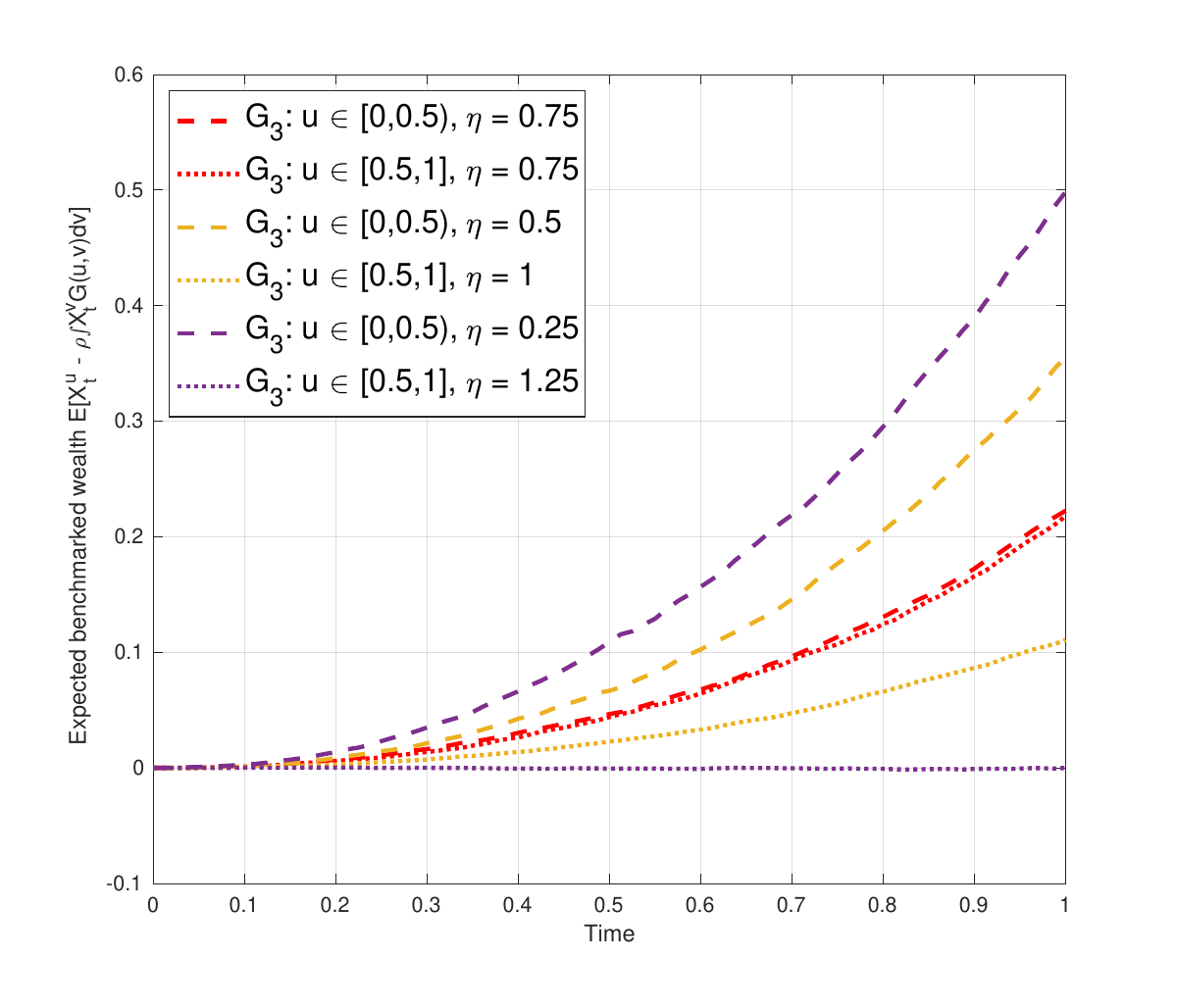}
\endminipage

\caption{Expected wealth and benchmarked wealth over time for $G_3$ when $\eta^u$ is nonconstant}
\label{Figure: expected wealth and relative wealth - G3 diff eta}
\end{figure}

For the minmax graphon, we plotted the expected wealth and expected benchmark wealth for $\eta^u = 1$ in red below as a benchmark for comparison. We experimented with taking $\eta^u$ to be a function $f(u) := \beta u(1-u)$, where $\beta$ indicates a scaling factor. Here we plotted two separate lines for the average wealth and average benchmarked wealth over time, for players with $u \in [0.25, 0.75]$ and players with $u \in [0,0.25) \cup (0.75,1]$. It is worth noting that unlike the star graphon case where the players are indistinguishable within each group, for the min-max graphon, the average wealth within each group, are just for comparison purposes for players with different labels. They are not representative of the expected wealth of any single player with a specific label. Considering the graphon $G_4$, and that $u(1-u)$ is concave and symmetric around $u = 0.5$, players with labels closer to $0.5$ will face a higher value for interaction strength $G(u,v)$, and a higher risk aversion parameter. As illustrated in \Cref{Figure: expected wealth and relative wealth - G4 diff eta} below, players with $u$ closer to $0.5$ have lower expected wealth and expected benchmark wealth. Meanwhile, higher values for the scaling factor $\beta$ lead to higher expected wealth and expected benchmark wealth for all players. We conducted similar experiments for the power-law graphon $G_5$ by taking $\eta^u = \beta u$, and plotted the average the plots are shown in \Cref{Figure: expected wealth and relative wealth - G5 diff eta} below. The effect of the scaling factor $\beta$ is consistant across different graphons.
\begin{figure}[H]
\minipage{0.49\textwidth}
    \includegraphics[width = \linewidth, right]{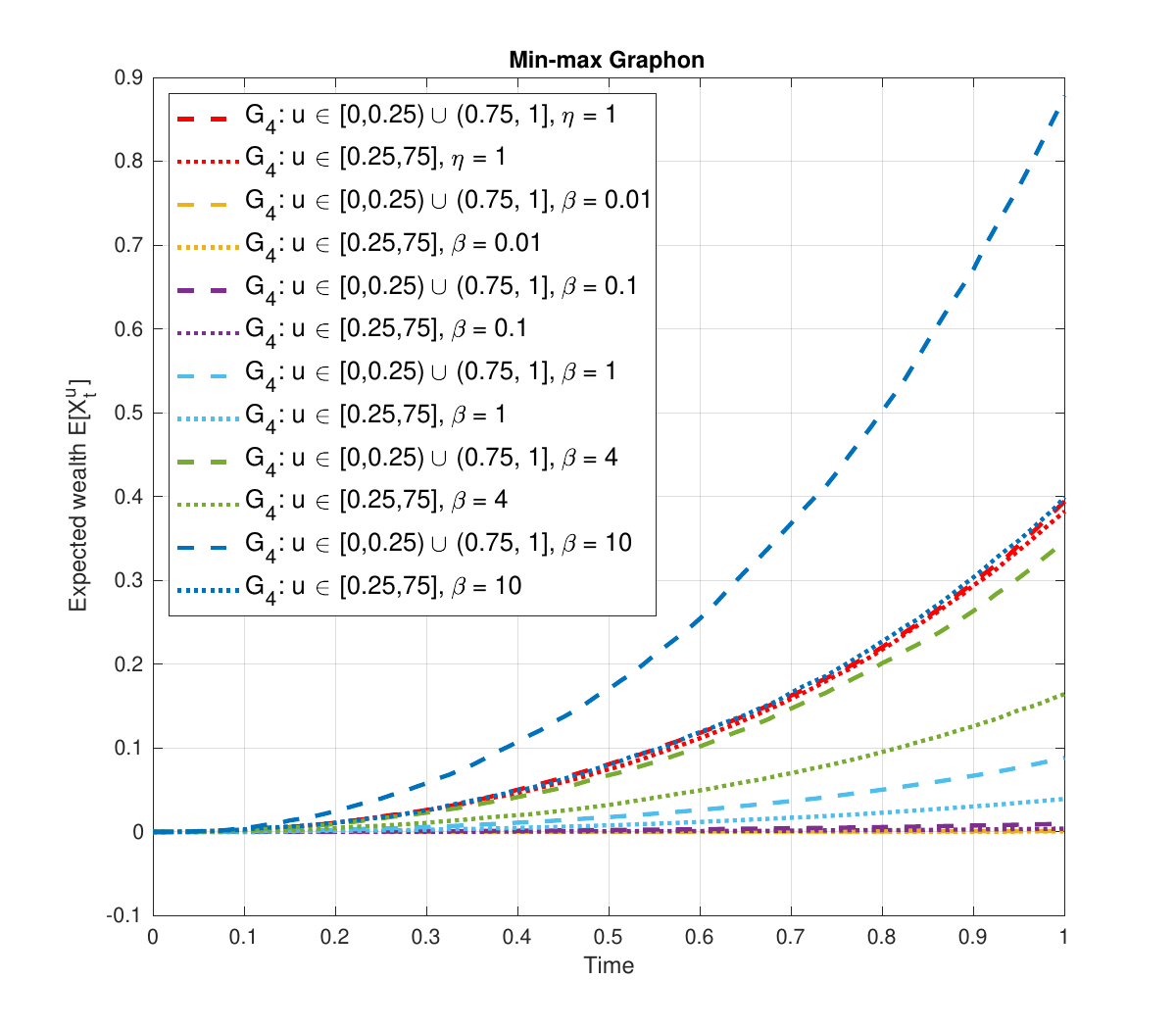}
\endminipage
\minipage{0.49\textwidth}
    \includegraphics[width = \linewidth, left]{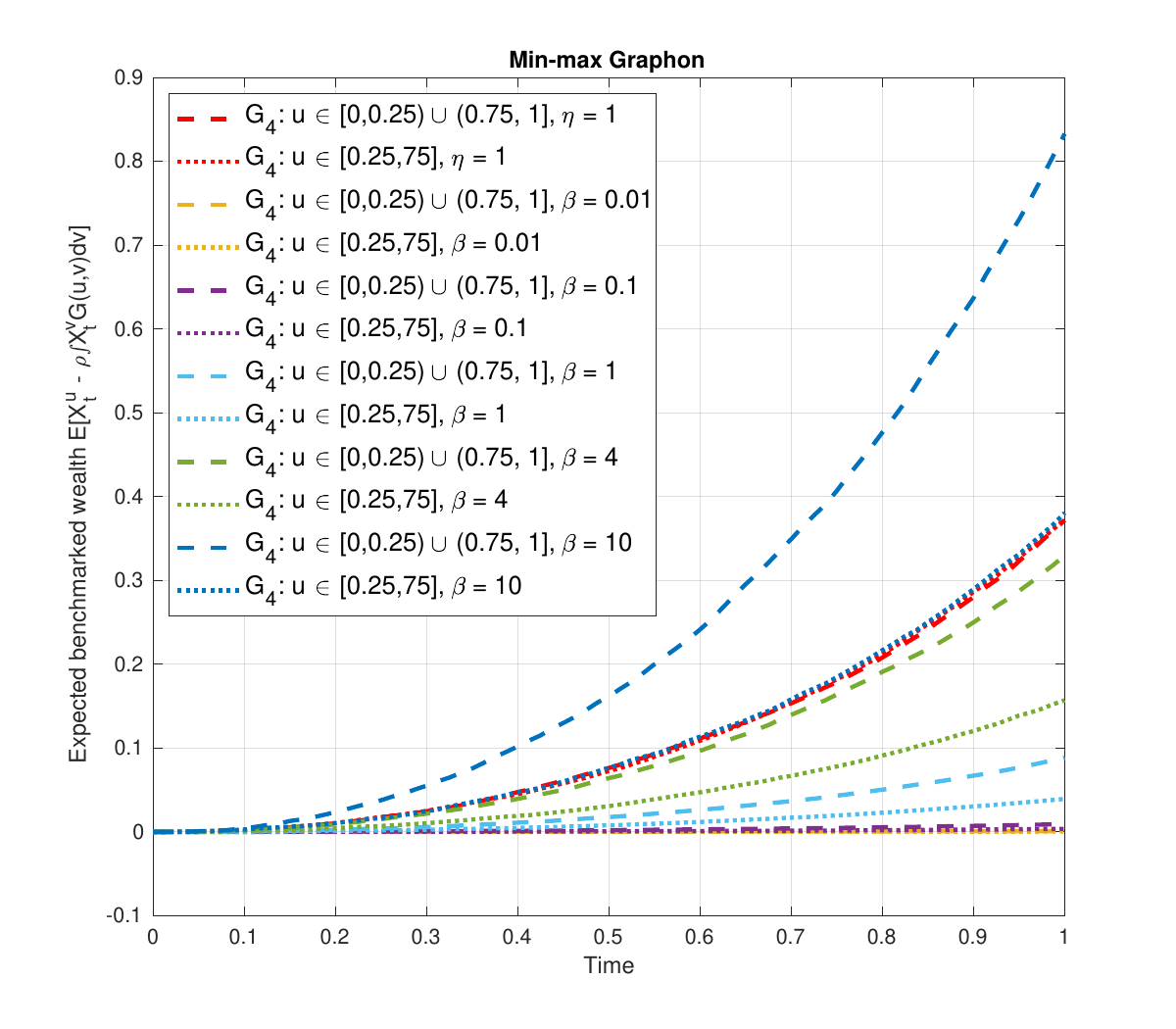}
\endminipage

\caption{Expected wealth and benchmarked wealth over time for $G_4$ when $\eta^u = \beta u(1-u)$}
\label{Figure: expected wealth and relative wealth - G4 diff eta}
\end{figure}

\begin{figure}[H]
\minipage{0.49\textwidth}
    \includegraphics[width = \linewidth, right]{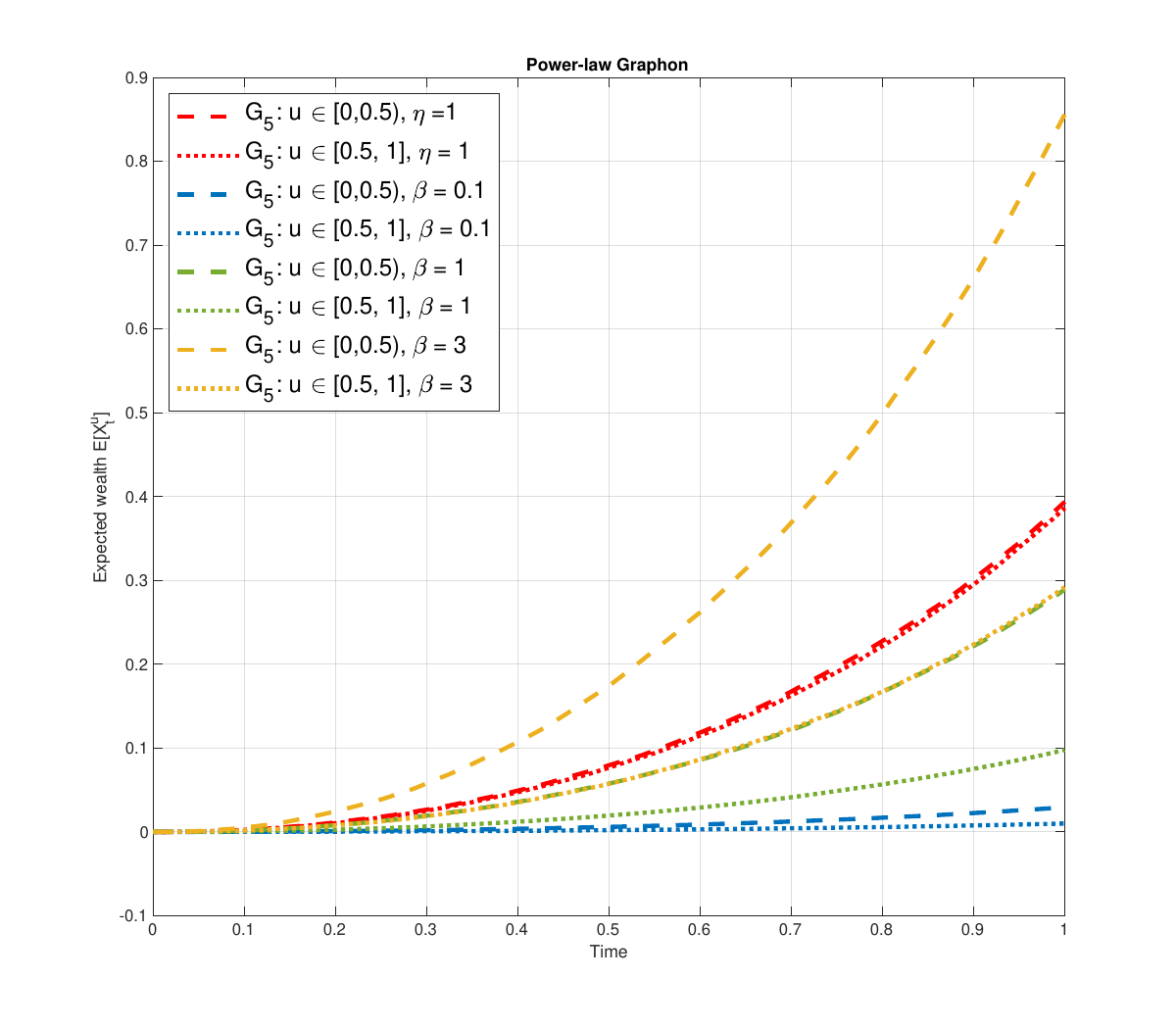}
\endminipage
\minipage{0.49\textwidth}
    \includegraphics[width = \linewidth, left]{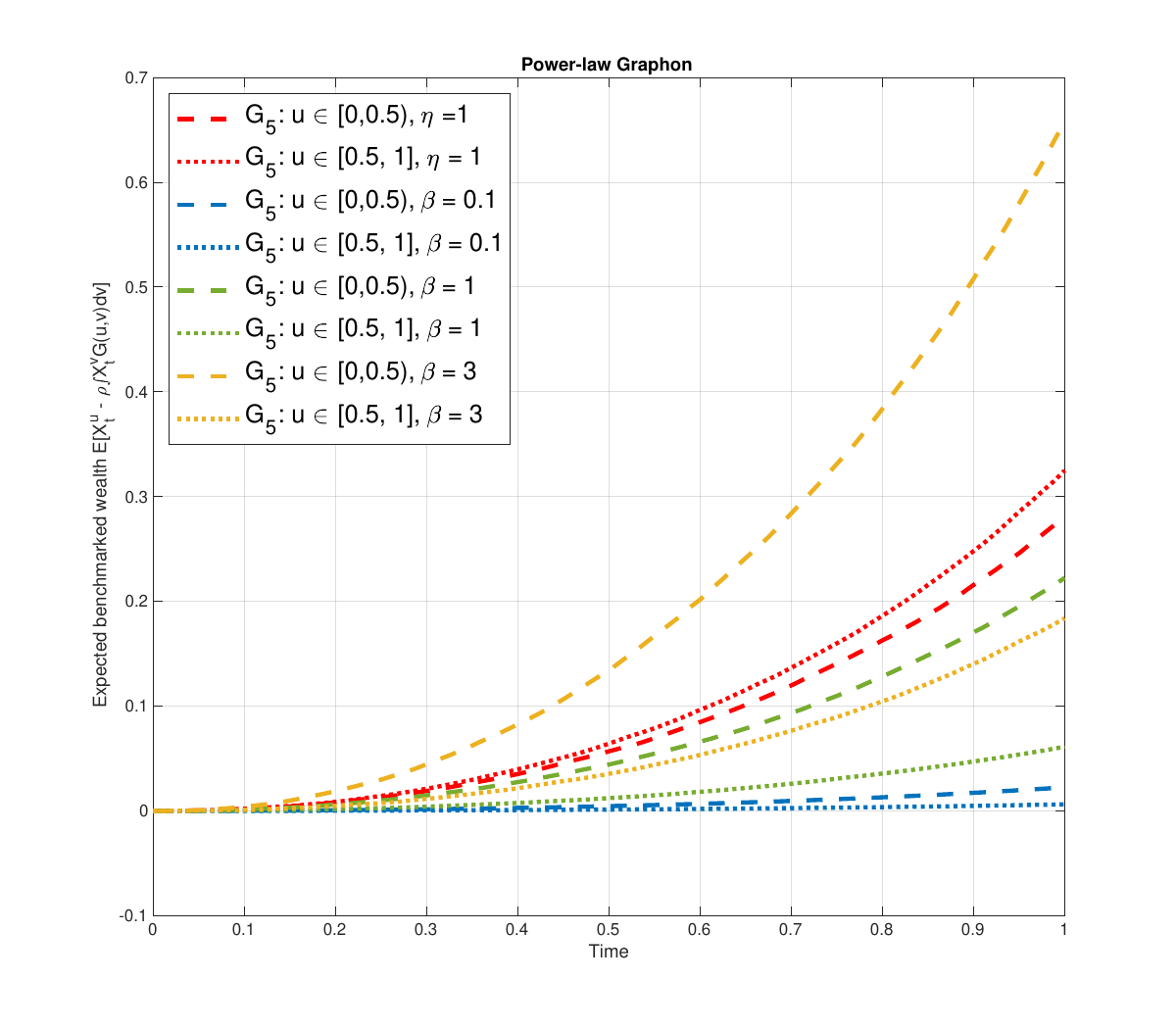}
\endminipage

\caption{Expected wealth and benchmarked wealth over time for $G_5$ when $\eta^u = u$}
\label{Figure: expected wealth and relative wealth - G5 diff eta}
\end{figure}

We would also like to investigate the effect of the risk-aversion parameter on the equilibrium utility \eqref{eq:graphon-obj}. 
We plotted the equilibrium utilities vs. labels for the min-max graphon and the power-law graphon. For the min-max graphon, the graph indicates that when $\eta^u$ is constant, $V_0^{u,G}$ is a convex function of $u$. An intuitive explanation is that players with labels closer to $0.5$ face higher interaction strength, hence more competition pressure, resulting in lower utilities. Notice that by taking $\eta^u = \beta u(1-u)$ and gradually increasing the scaling factor $\beta$, we first observe a change in concavity of $V_0^{u,G}$. When $\beta$ becomes large enough, the effect of the graphon on the relationship between utility and labels becomes negligible. Meanwhile, despite the fact that the expected terminal benchmarked wealth increases when $\beta$ increases, as illustrated in \Cref{Figure: expected wealth and relative wealth - G4 diff eta}, the equilibrium utilities remain unchanged, as illustrated below where the utility curves coincide when $\beta = 1$, $4$ and $10$. 
The plot for the power-law graphon shows similar results. When $\eta^u$ is constant, utilities are lower for higher values of $u$ since the interaction strength is stronger for players with these labels. When we take $\eta^u = \beta u$, as the scaling factor $\beta$ grows, players with higher values of $u$ gets higher utilities, and the utility curve as a function of $u$ remains unchanged after $\beta$ gets large enough. 

In conclusion, our experiment with varying the risk-aversion parameter indicates that players that are more risk-averse in general have higher expected terminal wealth and higher expected utility up to a certain extent. When $\eta^u$ becomes large enough, the expected utilities remain the same with further increase in scale for the risk-aversion parameter, despite the fact that the expected terminal benchmarked wealth kept increasing.
\begin{figure}[H]
\minipage{0.49\textwidth}
    \includegraphics[width = \linewidth, right]{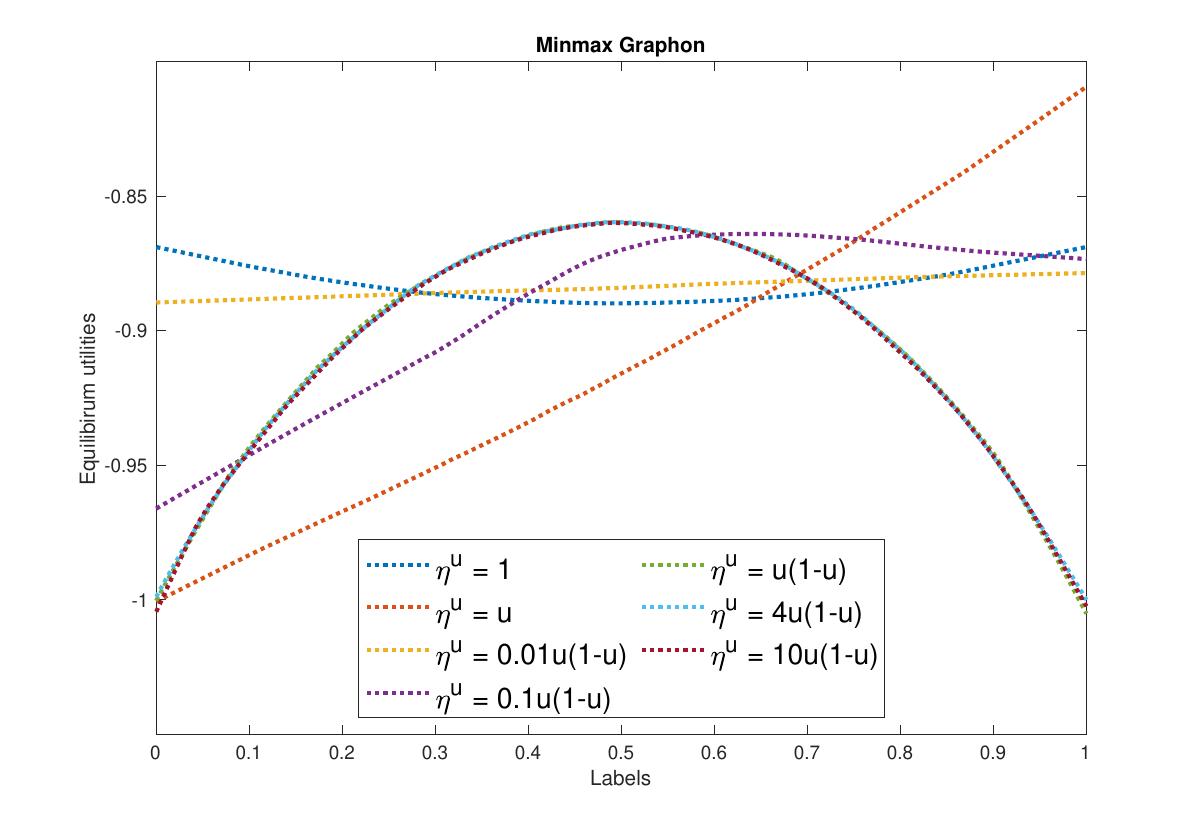}
\endminipage
\minipage{0.49\textwidth}
    \includegraphics[width = \linewidth, left]{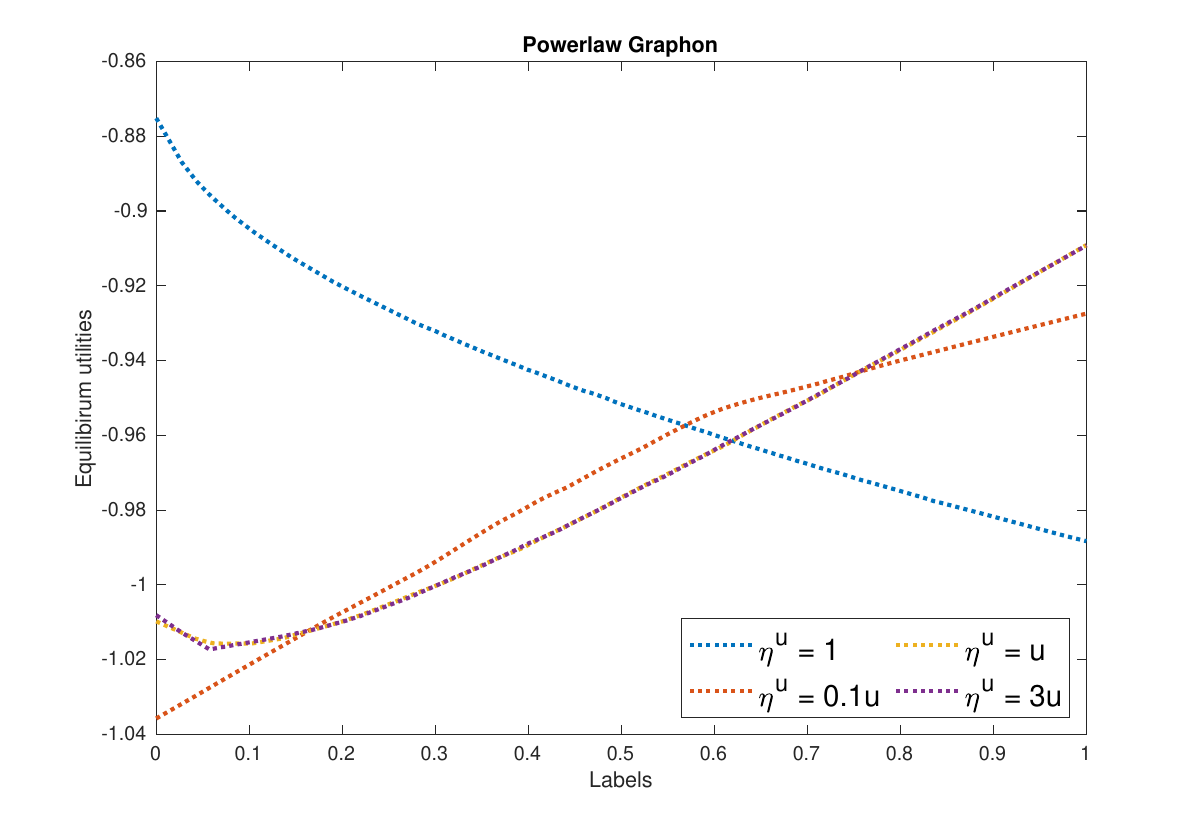}
\endminipage

\caption{Equilibrium utilities for $G_4$ and $G_5$ when $\eta^u$ is nonconstant in $u$}
\label{Figure: utility vs. label fn eta}
\end{figure}

\bibliographystyle{abbrvnat}
\bibliography{ref}

\section{Appendix}\label{sec: appendix}
As announced in Section \ref{sec:BS.Markov}, here we present the simulation result for the path of the value process $(Y_t)_{t\in [0,T]}$ for the Black-Scholes model with random coefficients.
\begin{figure}[H]
\minipage{0.33\textwidth}
    \includegraphics[width = \linewidth]{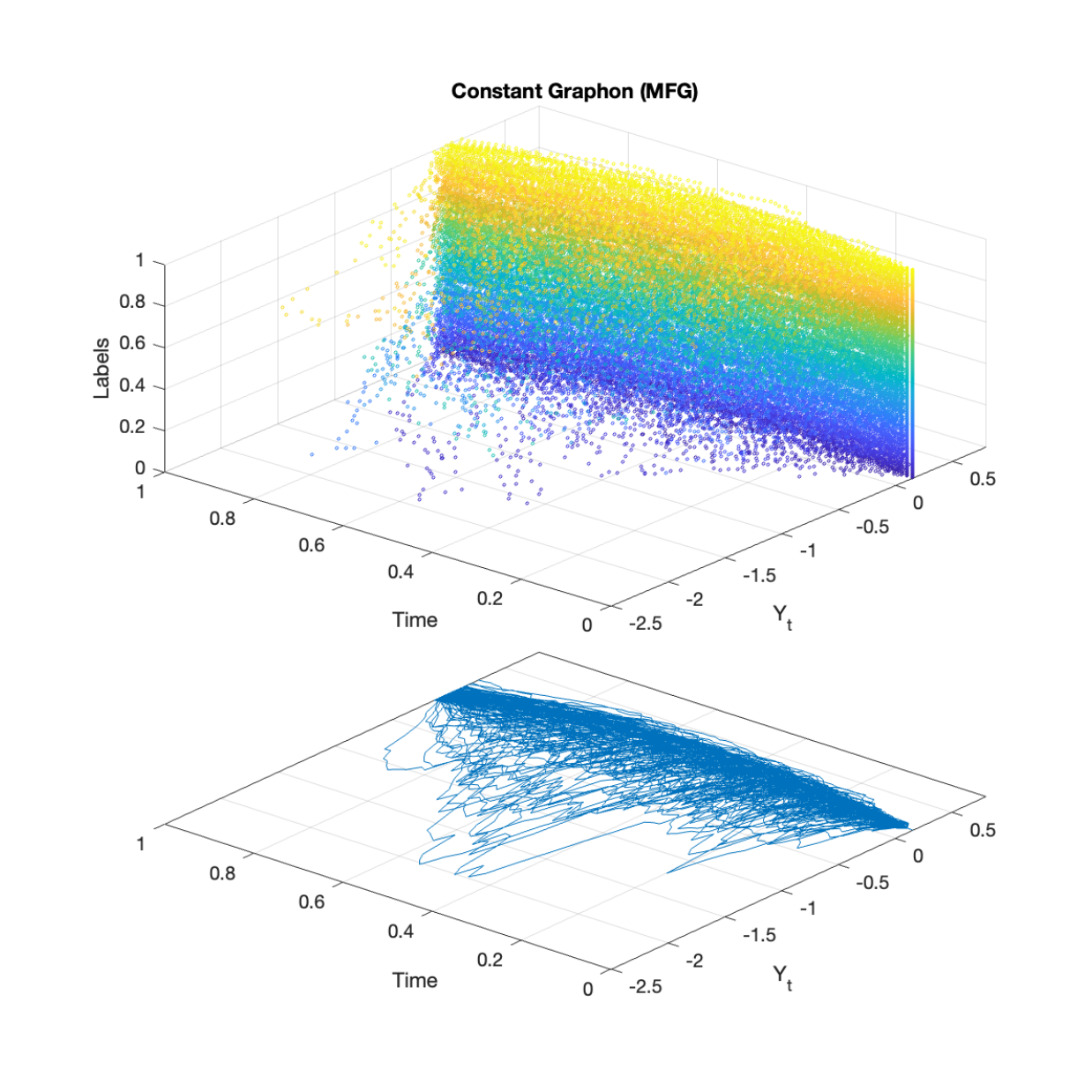}
\endminipage
\minipage{0.33\textwidth}
    \includegraphics[width = \linewidth]{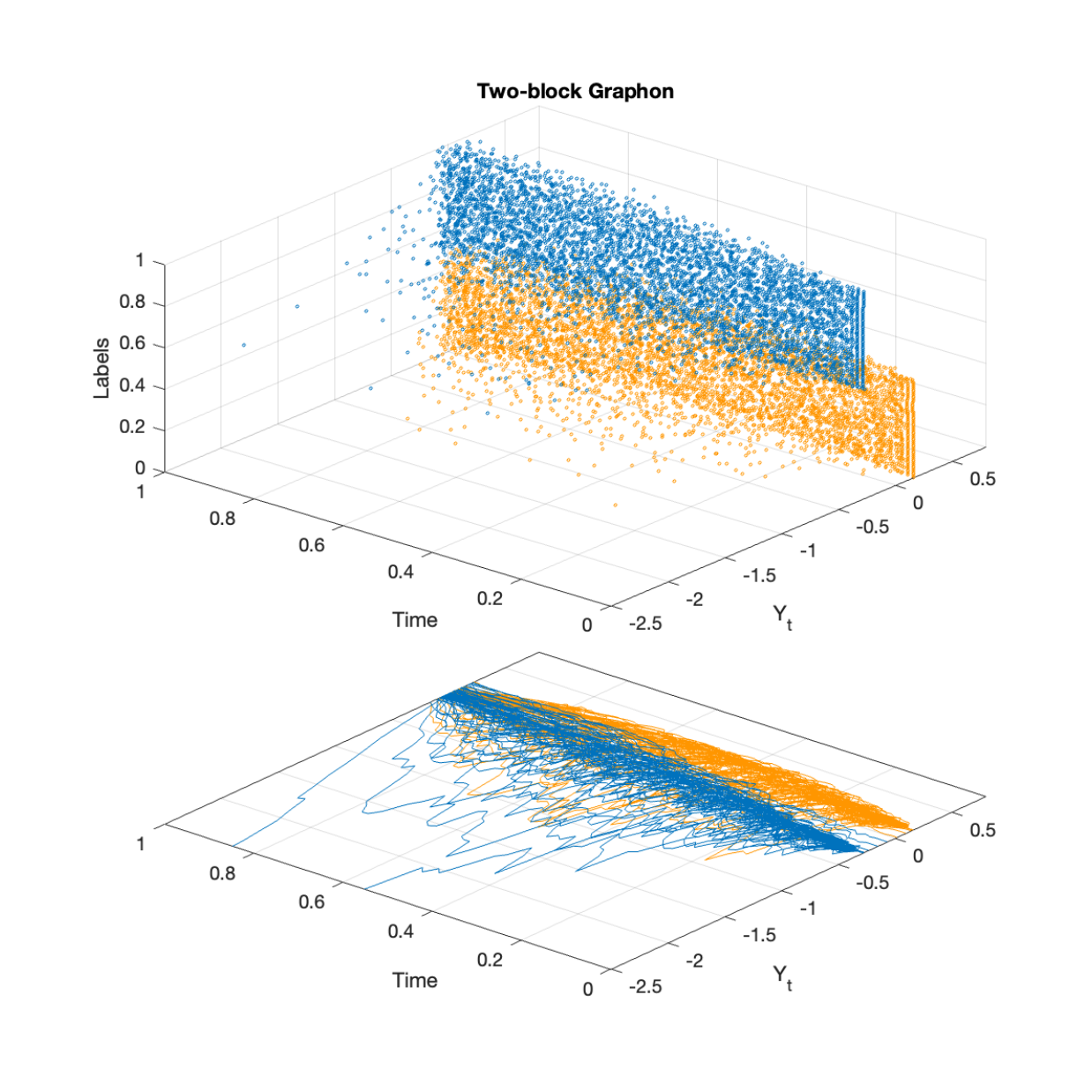}
\endminipage
\minipage{0.33\textwidth}
    \includegraphics[width = \linewidth]{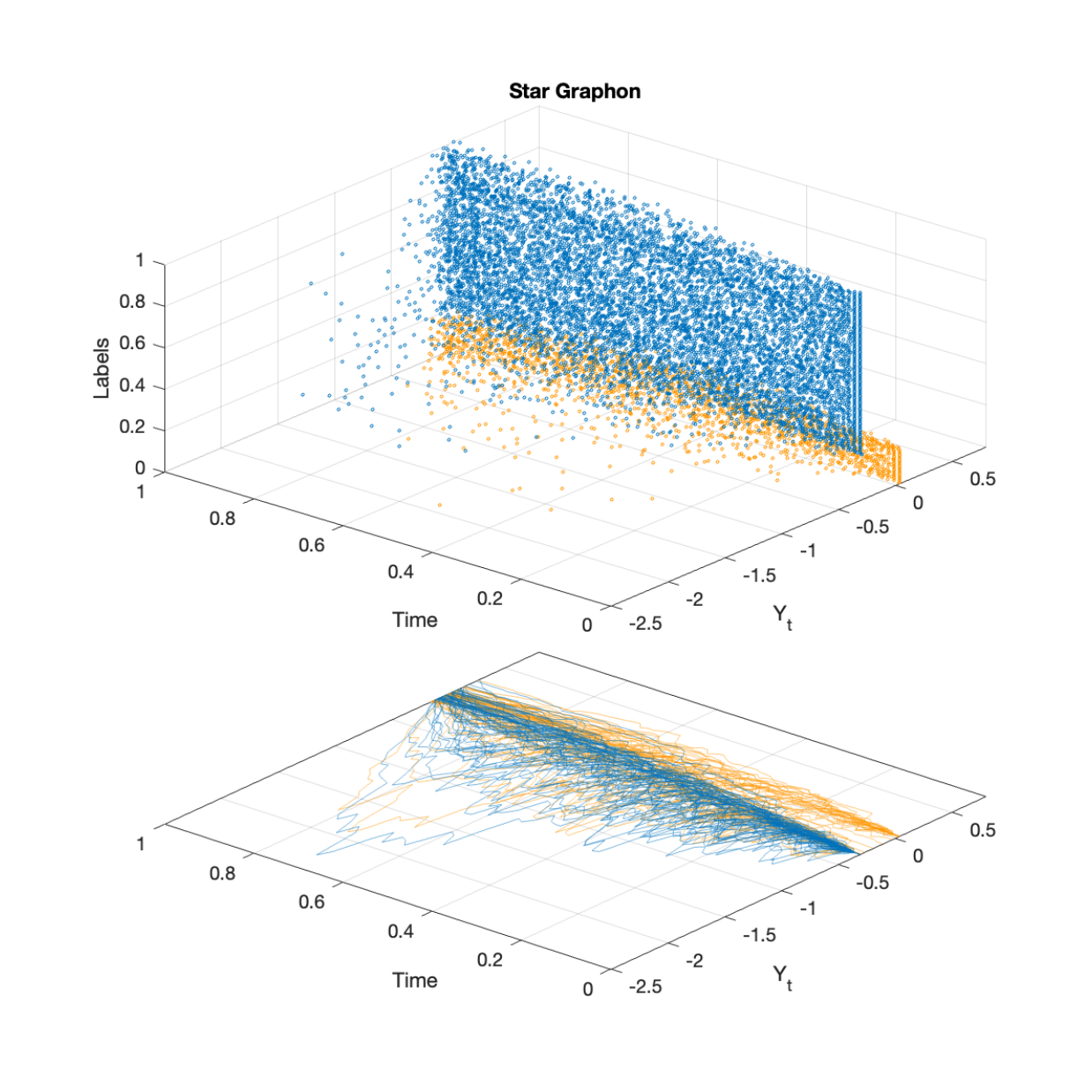}
\endminipage

\minipage{0.49\textwidth}
    \includegraphics[width = .75\linewidth, right]{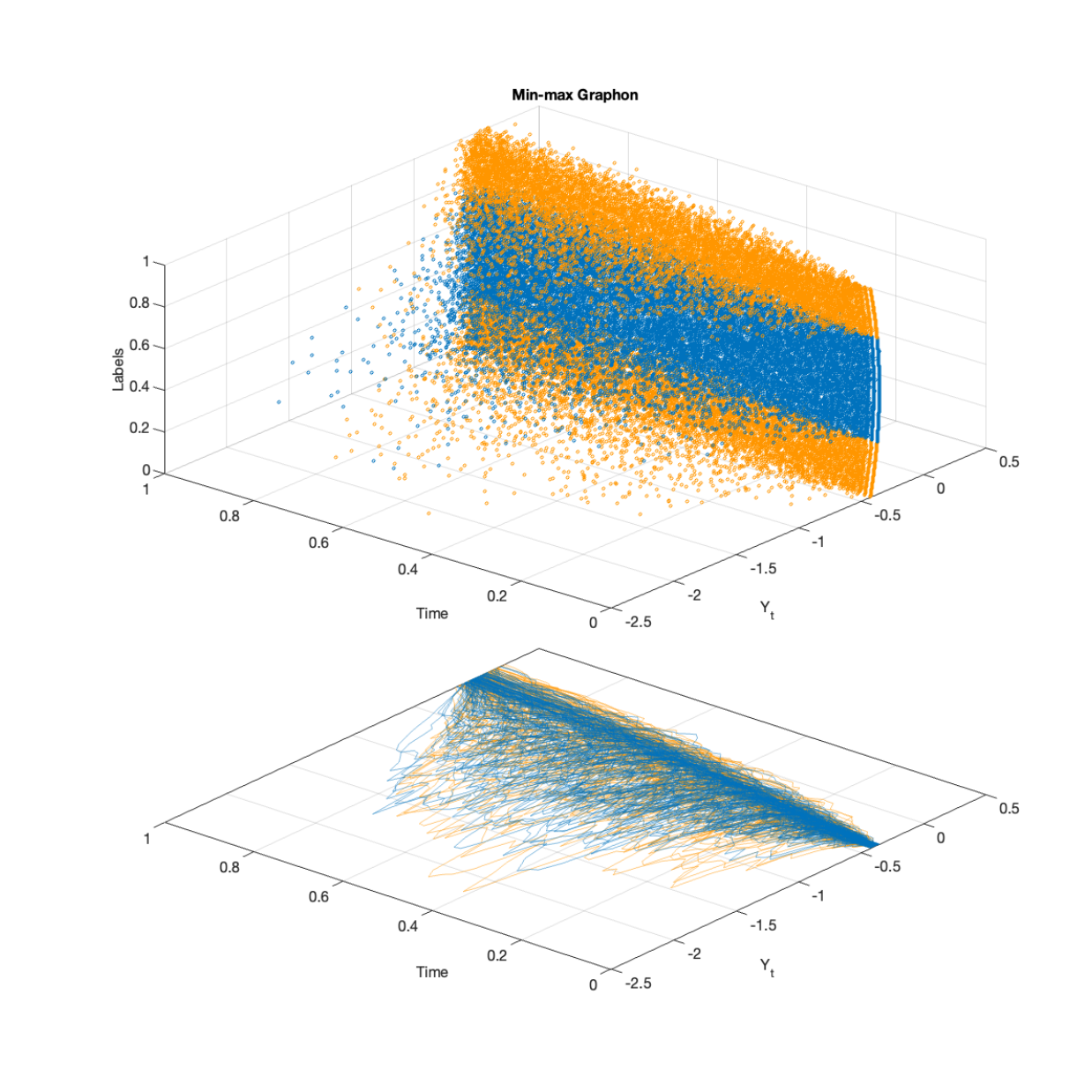}
\endminipage
\minipage{0.49\textwidth}
    \includegraphics[width = .75\linewidth, left]{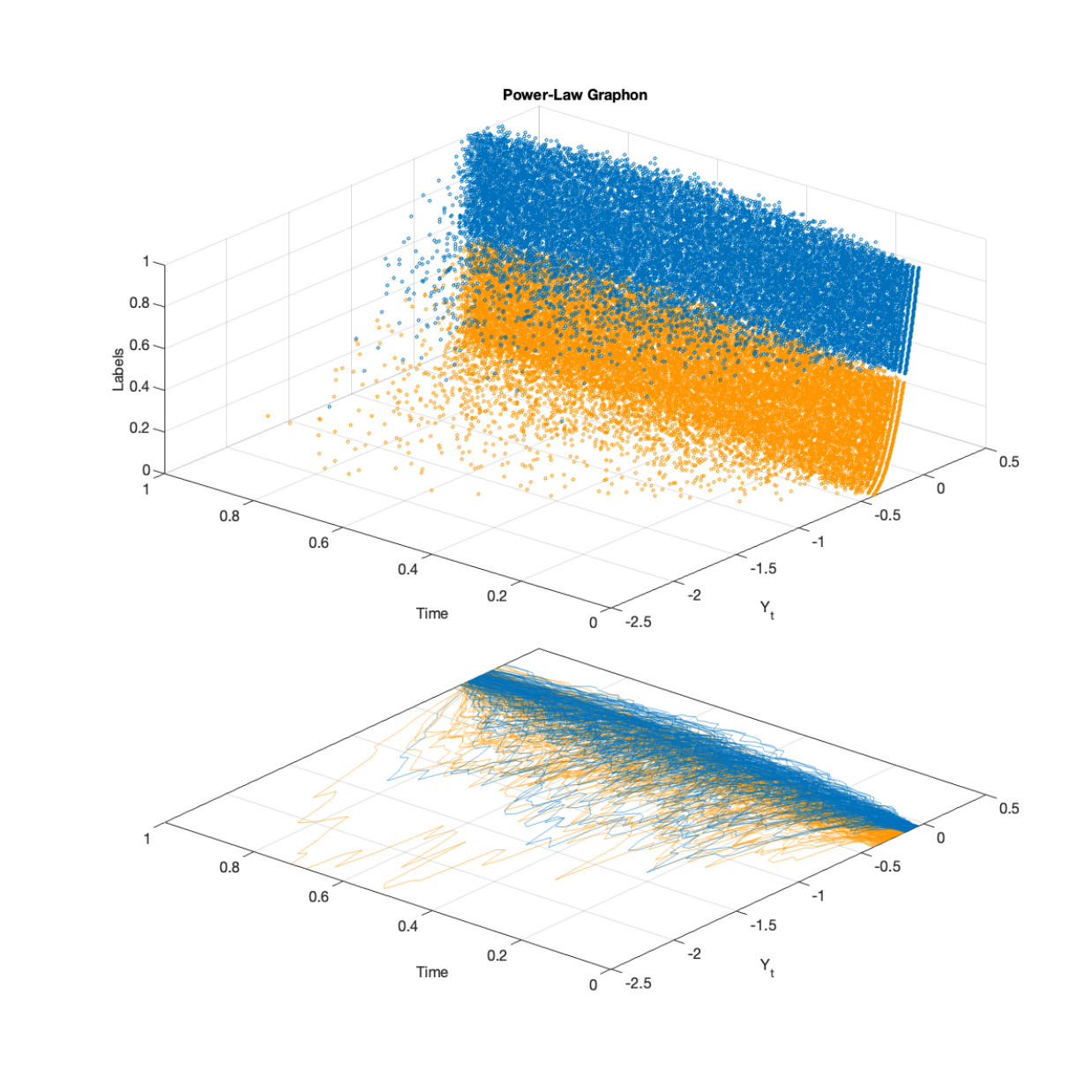}
\endminipage
\caption{Simulation of path $Y_t$ in graphon equilibrium}
\label{Figure: markovian BS}
\end{figure}

\end{document}